\RequirePackage{fix-cm}
\documentclass[smallextended]{svjour3}       
\smartqed  
\usepackage{graphicx}
\usepackage{dsfont}
\usepackage{amsmath}
\usepackage{float}
\usepackage{subfigure}
\usepackage{multirow}
\usepackage{multicol}
\usepackage[figuresright]{rotating}
\usepackage{epstopdf}
\newtheorem{assumption}{Assumption}

\journalname{Computational Optimization and Applications}
\begin{document}

\title{Mini-batch stochastic Nesterov's smoothing method for constrained convex stochastic composite optimization
\thanks{The work is supported in part by ``the Fundamental Research Funds for the Central Universities" (Grant No. 2019YJS199) and ``the Natural Science Foundation of Beijing, China" (Grant No. 1202021)}}

\titlerunning{Mini-batch stochastic Nesterov's smoothing method}        

\author{Ruyu Wang         \and
        Chao Zhang        \and
        Lichun Wang \and
        Yuanhai Shao
}


\institute{Ruyu Wang \at
           Department of Applied Mathematics, Beijing Jiaotong University \\
              \email{wangruyu@bjtu.edu.cn}           
           \and
           Chao Zhang \at
           Department of Applied Mathematics, Beijing Jiaotong University \\
              \email{zc.njtu@163.com}           
           \and
           Lichun Wang \at
           Department of Applied Mathematics, Beijing Jiaotong University \\
              \email{lchwang@bjtu.edu.cn}
              \and
           Yuanhai Shao \at
           Management School, Hainan University \\
              \email{shaoyuanhai@hainanu.edu.cn}
}

\date{Received: date / Accepted: date}

\maketitle

\begin{abstract}
This paper considers a class of constrained convex stochastic composite optimization problems whose objective function is given by the summation of a differentiable convex component, together with a nonsmooth but convex component. The nonsmooth component has an explicit max structure that may not easy to compute its proximal mapping. In order to solve these problems, we propose a mini-batch stochastic Nesterov's smoothing (MSNS) method. Convergence and the optimal iteration complexity of the method are established. Numerical results are provided to illustrate the efficiency of the proposed MSNS method for a support vector machine (SVM) model.

\keywords{Constrained convex stochastic programming\and Mini-batch of samples\and Stochastic approximation\and Nesterov's smoothing method\and Complexity}
\end{abstract}

\section{Introduction}
In this paper, we consider the nonsmooth convex stochastic composite minimization problem
\begin{equation}\label{orip}
\psi^*:=\min\limits_{x\in X}\{\psi(x):=f(x)+h(x)\},
\end{equation}
where $X$ is a bounded closed convex  set in the Euclidean space $\mathds{R}^n$, $f:X\rightarrow\mathds{R}$ is a convex function with Lipschitz continuous gradient, and $h:X\rightarrow\mathds{R}$ is a possibly nonsmooth convex function with the explicit max structure
\begin{eqnarray}\label{h}
h(x) = \max_{u\in U} \{\langle Ax, u\rangle - Q(u) \},
\end{eqnarray}
where $U$ is a bounded closed convex set, $A$ is a linear operator, and $Q(u)$ is a continuous convex function. Such function $h$ has been studied by Nesterov \cite{Nesterov} with various important applications, and its max structure has been used to construct the smoothing approximation $h_{\mu}(x)$ of $h(x)$ so that its gradient is Lipschitz continuous.

Although $f$ and $h_{\mu}$ are Lipschitz continuously differentiable, we assume that only the noisy objective values and gradients of $f$ and $h_{\mu}$ are available via subsequent calls to a stochastic oracle ($\cal SO$). Many applications especially in machine learning are in this setting. That is, when we solve the smooth problem
\begin{eqnarray}\label{smooth}
\min\limits_{x\in X}\{\psi_{\mu}(x):=f(x)+h_{\mu}(x)\},
\end{eqnarray}
by an iterative algorithm, at the $k$-th iterate, $k\ge 1$, for the input $x_k\in X$, the $\cal SO$ would output a stochastic value $\Psi_{\mu}(x_k,\xi_k)$ and a stochastic gradient $\nabla\Psi_{\mu}(x_k,\xi_k)$ in the form of
\begin{eqnarray*}
\Psi_{\mu}(x_k, \xi_k):=F(x_k, \xi_k)+H_{\mu}(x_k, \xi_k)
\end{eqnarray*}
and
\begin{eqnarray*}
\nabla \Psi_{\mu}(x_k, \xi_k):=\nabla F(x_k, \xi_k)+\nabla H_{\mu}(x_k, \xi_k),
\end{eqnarray*}
where $\xi_k$ is a random vector whose probability distribution is supported on $\Xi\subseteq\mathds{R}^d$.

\begin{assumption}\label{assumption 1}
For any fixed $\mu>0$, and $x\in X$, we have
\begin{eqnarray*}
&a)&E[\Psi_{\mu}(x, \xi)]=\psi_{\mu}(x),\\
&b)&E[\nabla\Psi_{\mu}(x, \xi)]=\nabla\psi_{\mu}(x),\\
&c)&E[\|\nabla\Psi_{\mu}(x, \xi)-\nabla\psi_{\mu}(x)\|^2]\leq \sigma^2,
\end{eqnarray*}
where $\sigma>0$ is a constant, and the expectation $E$ is taken with respect to the random vector $\xi\in \Xi$.
\end{assumption}

The stochastic approximation (SA)  is one important approach for solving stochastic convex programming, which can be dated back to the pioneering paper by Robbins and Monro \cite{Robbins}. A robust version of the SA method developed by Polyak \cite{Polyak}, and Polyak and Juditsky \cite{Polyak2}, improves the original version of the SA method. It was demonstrated in Nemirovski et al. \cite{Nemirovski1} that a proper modification of the SA approach based on the mirror-descent SA (Nemirovski and Yudin \cite{Nemirovski2}), can be competitive and can even significantly outperform the other important type approach, the sample average approximation (SAA) method  \cite{Kleywegt,Shapiro1}, for a certain class of convex stochastic programming in \cite{Lan,Nemirovski}. It has been pointed out in \cite{Nemirovski1} that for nonsmooth stochastic convex optimization, the iteration complexity of order
$$
O\left(\frac{1}{\epsilon^2}\right)
$$
is optimal, where $\epsilon$ is the desired absolute accuracy of the approximate solution in objective value.
In 2016, S. Ghadimi et al. \cite{Ghadimi} proposed a novel randomized stochastic projected gradient (RSPG) algorithm which can solve constrained nonconvex nonsmooth stochastic composite problems. 
The nonsmooth component $h$ in \cite{Ghadimi}, however, is restricted to be a simple convex function such as $h(x)=\|x\|_1$, for which the proximal operator is easy to compute. Such restriction still exists for the other stochastic proximal type methods, e.g., \cite{Wang,Ma,Xiao}.
There are many real applications in which the nonsmooth term is not so easy to obtain the proximal operator, such as  a support vector machine (SVM) model we consider in numerical experiment in Section 4, where $h$ is a maximum of $0$ and an affine function.

Smoothing methods have been shown to be efficient for dealing with constrained nonsmooth optimization with solid convergence results, which allow the nonsmooth terms to be relatively complex \cite{Nesterov,Zhang,Chen,Bian,Liu,Zhang2}. In this paper, we propose a mini-batch stochastic Nesterov's smoothing (MSNS) method for solving \eqref{orip} with relatively complex nonsmooth convex component $h$. Note that Nesterov's smoothing method  \cite{Nesterov} was designed for solving deterministic constrained convex nonsmooth composite optimization problems. The MSNS method proposed in this paper is suitable for stochastic setting. The extension, however, is not a trivial task. We show the convergence, as well as the optimal iteration complexity of the MSNS method.  We illustrate the efficiency of the MSNS method, by comparing with several state-of-the-art SA-type methods, on a SVM model using both synthetic data and several real data.

The remaining part of this paper is organized as follows. In Section 2, we briefly review some basic concepts and results relating to the Nesterov's smoothing method \cite{Nesterov} that will be used in our paper. In Section 3, we develop a mini-batch stochastic Nesterov's smoothing method which extends the Nesterov's smoothing method from the deterministic setting to the stochastic setting. We show the convergence, as well as the optimal iteration complexity of the proposed method. Numerical experiments on a SVM model are given in Section 4 to demonstrate the efficiency of our proposed method.

\section{Preliminaries}

In this section, we review some basic concepts and results relating to the Nesterov's smoothing method \cite{Nesterov} that will be used later. The problem \eqref{orip} can be considered as a convex-concave saddle
point problem $$\min\limits_{x\in X}\max\limits_{u\in U} \left\{K(x,u):=f(x)+\langle Ax,u\rangle-Q(u)\right\},$$
where $K$ is convex in $x$ on $X$, and concave in $u$ on $U$.
The adjoint form of the problem \eqref{orip} can be written as
\begin{eqnarray}\label{eq2.5}
\max\limits_{u\in U}\left\{\phi(u):=\min\limits_{x\in X}K(x,u)\right\}.
\end{eqnarray}
We call \eqref{orip} the ``primal problem" and \eqref{eq2.5} the ``dual problem". It is easy to see that
\begin{eqnarray*}
\psi(x) \ge \phi(u), \quad \mbox{for any}\ x\in X\ \mbox{and}\ u\in U.
\end{eqnarray*}
Since the objective function of each problem is continuous and the feasible region is compact, we know that
there exists $x^*\in X$ and $u^*\in U$, which are optimal solutions of \eqref{orip} and \eqref{eq2.5}, respectively. According to Theorem 4.2' in \cite{Sion}, we have
\begin{eqnarray*}
\min\limits_{x\in X}\max\limits_{u\in U}K(x,u)=\max\limits_{u\in U}\min\limits_{x\in X}K(x,u).
\end{eqnarray*}
That is, $\psi(x^*) = \phi(u^*)$. Hence in the Nesterov's smoothing method \cite{Nesterov} for the determinstic setting,
if the gap $\psi({\widehat x})-\phi({\widehat u}) \le \epsilon$ for a tolerance $\epsilon>0$, the output $\widehat x$ is called an $\epsilon$-approximate solution of the primal problem \eqref{orip}.

As in \cite{Nesterov}, by inserting a non-negative, continuous and $\sigma_\omega$-strongly convex function $\omega(u)$ in \eqref{h}, we obtain a smooth approximation $h_{\mu}(x)$ of $h(x)$
\begin{eqnarray}\label{eq2.6}
h_{\mu}(x):=\max\limits_{u\in U}\left\{\langle Ax,u\rangle-Q(u)-\mu\omega(u)\right\},
\end{eqnarray}
where $\mu>0$ is a smoothing parameter. Let us denote by $u_{\mu}(x)$ the optimal solution of the above problem. Recall that  $\omega(u)$ is strongly convex on $U$ if there exists a constant $\sigma_\omega>0$ such that
\begin{eqnarray}\label{eq2.10}
\omega(u_1)\geq \omega(u_2)+\left\langle\nabla\omega(u_2),u_1-u_2 \right\rangle+ \frac{1}{2}\sigma_\omega\|u_1-u_2\|^2, ~~\forall u_1,u_2\in U.
\end{eqnarray}

Denote by $u_0=\arg\min\limits_{u\in U}\{\omega(u)\}$. Without loss of generality we assume that $\omega(u_0)=0$. By \eqref{eq2.10}, for any $u\in U$ we have
\begin{eqnarray}\label{eq2.8}
\omega(u)\geq\frac{1}{2}\sigma_\omega\|u-u_0\|^2.
\end{eqnarray}
Let $\Omega=\max\limits_{u\in U}\{\omega(u)\}$. Then, according to \cite{Nesterov}, for any $x\in X$ we have
\begin{eqnarray}\label{eq2.7}
h_{\mu}(x)\leq h(x)\leq h_{\mu}(x)+\mu \Omega.
\end{eqnarray}

\begin{lemma}\label{lemm2.1}(Theorem 1 of \cite{Nesterov})
The function $h_{\mu}(x)$ is well defined and continuously differentiable at any $x\in X$. Moreover, this function is convex and its gradient
\begin{eqnarray*}
\nabla h_{\mu}(x)=A^Tu_{\mu}(x)
\end{eqnarray*}
is Lipschitz continuous with constant
\begin{eqnarray}\label{eq2.9}
L_{h_\mu}=\frac{1}{\mu \sigma_\omega}\|A\|^2,
\end{eqnarray}
where $\|A\|$ is the operator norm of $A$.
\end{lemma}
We then find $\psi_{\mu}(x)$ is $L$-smooth on $X$, i.e., $\psi_{\mu}(x)$ is Lipschitz continuously differentiable with Lipschitz constant $0<L=L_f+L_{h_{\mu}}$ on $X$. According to Theorem 5.12 of \cite{Beck}, we get $L_f\geq\|\nabla^2f(x)\|=\lambda_{\max}(\nabla^2f(x))$ for any $x\in X$, where $\lambda_{\max}(\nabla^2f(x))$ means the maximal eigenvalue of the Hessian matrix of $f(x)$. Thus
\begin{eqnarray}\label{eq3.10}
\psi_{\mu}(y)\leq \psi_{\mu}(x)+\langle \nabla \psi_{\mu}(x),y-x\rangle+\frac{L}{2}\|y-x\|^2,\quad \forall x, y \in X.
\end{eqnarray}

Let $d(x)$ be a prox-function of the set $X$. We assume that $d(x)$ is continuous and
strongly convex on $X$ with convexity parameter $\sigma_d > 0$. Let $x_0$ be the center of the set $X$, i.e.,
$x_0=\arg\min \{d(x) : x\in X\}$. Without loss of generality we assume that $d(x_0)=0$. Thus
\begin{eqnarray}\label{eq3.11}
d(x) \geq \frac{1}{2}\sigma_d\|x-x_0\|^2.
\end{eqnarray}
In view of the compact set $X$, we know that there exists a constant $D>0$ such that
\begin{eqnarray}\label{eq3.12}
\max\limits_{x\in X}\{d(x)\} \leq D.
\end{eqnarray}
Given $g\in \mathds{R}^n$, we define the generalized projected gradient of $\psi$ at $x \in X$ as
\begin{eqnarray}\label{eq3.13}
P_X(x,g,\gamma)=\frac{1}{\gamma}(x-x'),
\end{eqnarray}
where $x'$ is given by
\begin{eqnarray}\label{x'}
x'=\arg\min\limits_{y\in X}\left\{\langle g,y\rangle+\frac{1}{2\gamma}\|y-x\|^2\right\}.
\end{eqnarray}

\begin{lemma}\label{Lemma2}
Let $P_X(x,g,\gamma)$ be defined in \eqref{eq3.13}-\eqref{x'}. Then, for any $x\in X$, $\mu>0$ and $\gamma>0$,
the stochastic gradient $\nabla\Psi_{\mu}(x,\xi)$ and the full gradient $\nabla\psi_{\mu}(x)$  satisfy
\begin{eqnarray}\label{Pleqg}
\|P_X(x,\nabla\Psi_{\mu}(x,\xi),\gamma)-P_X(x,\nabla\psi_{\mu}(x),\gamma)\|\leq \|\nabla\Psi_{\mu}(x,\xi)-\nabla\psi_{\mu}(x)\|.
\end{eqnarray}
\end{lemma}
\begin{proof}
The generalized projected gradient defined in \eqref{eq3.13}-\eqref{x'} is a special case of that defined in (2.3)-(2.4) of \cite{Ghadimi}, with $V(y,x) = 1/2\|y-x\|^2$, the convexity parameter $\alpha = 1$, and $h(x) = 0$ in (2.3) of \cite{Ghadimi}. Then the statement in this lemma can be obtained easily from Proposition 1 of \cite{Ghadimi}, which claims that the generalized projected gradient $P_X(x,\cdot,\gamma)$ is Lipschitz continuous with Lipschitz constant $1/{\alpha} = 1$.
\end{proof}

\section{Mini-batch Stochastic Nesterov's Smoothing Method}

In this section, we develop a novel mini-batch stochastic Nesterov's smoothing (MSNS) method for solving \eqref{orip}. We will show the convergence as well as the optimal iteration complexity of the MSNS method.

At the $k$-th iterate of the MSNS method, we randomly choose a mini-bath samples $\xi^{[m_k]}:=\{\xi_{k,1},\ldots,\xi_{k,m_k}\}$ of the random vector $\xi\in \Xi$, where $m_k$ is the batch size. And we denote by $\xi_{[k]}:=\{\xi^{[m_0]},\ldots,\xi^{[m_{k}]}\}$ the history of mini-batch samples from the $0$-th iterate up to the $k$-th iterate. For any $k\geq 0$, we denote the mini-batch stochastic objective value $\Psi_{\mu}^k$, and the mini-batch stochastic gradient $\nabla \Psi_{\mu}^k$ by
\begin{eqnarray}
\Psi_{\mu}^k = \frac{1}{m_{k}}\sum\limits_{j=1}^{m_{k}} \Psi_{\mu}(x_k, \xi_{k,j}), \quad   \nabla \Psi_{\mu}^k=\frac{1}{m_{k}}\sum\limits_{j=1}^{m_{k}}\nabla \Psi_{\mu}(x_k, \xi_{k,j}). \label{eq3.9}
\end{eqnarray}

For ease of notations in the proof, we denote
\begin{eqnarray}
& &\delta_{\mu}(x_k): = \nabla\Psi_{\mu}^k-\nabla\psi_{\mu}(x_k),\label{notation2-1}\\
& & \delta_{\mu}^i(x_k):= \nabla\Psi_{\mu}(x_k,\xi_{k,i})-\nabla\psi_{\mu}(x_k),\ i=1,\cdots,m_k, \label{notation2-2}\\
& & S_{\mu}^j(x_{k}) := \sum\limits_{i=1}^{j}\delta_{\mu}^i(x_k),\ j =1,\cdots,m_k. \label{notation2-3}
\end{eqnarray}

The following lemma addresses the relations of  $\Psi_{\mu}^k$ and $\nabla \Psi_{\mu}^k$ to that of the original objective value $\psi_\mu(x_k)$ and $\nabla \psi_{\mu}(x_k)$, which can be shown without difficulty using the arguments similar as in \cite{Ghadimi}.
\begin{lemma}\label{lemma3.3}
Under Assumption \ref{assumption 1}, we have for any $k\geq0$ and $\mu >0$,
\begin{eqnarray*}
&  (a)& E\left[\Psi_{\mu}^k\right] = E\left[\psi_{\mu}(x_k)\right],\\
&  (b)& E\left[\nabla \Psi_{\mu}^k\right] = E[\nabla \psi_{\mu}(x_k)] ,\\
&  (c)& E\left[\| \nabla \Psi_{\mu}^k - \psi_{\mu}(x_k) \|^2\right] \leq \frac{\sigma^2}{m_k},
\end{eqnarray*}
where the expectation is taken with respect to the history of mini-batch samples $\xi_{[k]}$.
\end{lemma}
\begin{proof}
Note that the $k$-th iterate $x_k$ is a function of the history $\xi_{[k-1]}$ of the generated random process and consequently is random, and it is independent of the mini-batch samples $\xi^{[m_k]}$ at the $k$-th iterate.

For $k=0$, no history of mini-batch samples  to call the $\cal SO$ exists before this iterate.
Thus Lemma \ref{lemma3.3} can be obtained immediately using Assumption 1.

Now we consider the case for any $k\ge 1$. Using Assumption \ref{assumption 1} (a), we have
\begin{eqnarray*}
E\left[\Psi_{\mu}^k\right]&=&E\left\{E\left[\Psi_{\mu}^k |\xi_{[k-1]}\right]\right\}\\
&\overset{\eqref{eq3.9}}{=}&E\left\{E\left[\frac{1}{m_{k}}\sum\limits_{j=1}^{m_{k}} \Psi_{\mu}(x_k, \xi_{k,j})| \xi_{[k-1]}\right]\right\}\\
&=&E\left\{\frac{1}{m_{k}}\sum\limits_{j=1}^{m_{k}} E\left[\Psi_{\mu}(x_k, \xi_{k,j}) |\xi_{[k-1]}\right]\right\}\\
&{=}&E\left[\psi_{\mu}(x_k)\right].
\end{eqnarray*}
Employing Assumption \ref{assumption 1} (b), we get
\begin{eqnarray*}
E\left[\nabla\Psi_{\mu}^k\right]
&=&E\left\{E\left[\nabla\Psi_{\mu}^k|\xi_{[k-1]}\right]\right\}\\
&\overset{\eqref{eq3.9}}{=}&E\left\{ E\left[\frac{1}{m_{k}}\sum\limits_{j=1}^{m_{k}}\nabla \Psi_{\mu}(x_k, \xi_{k,j}) | \xi_{[k-1]}\right] \right\}\\
&=& E\left\{\frac{1}{m_{k}}\sum\limits_{j=1}^{m_{k}}E\left[\nabla \Psi_{\mu}(x_k, \xi_{k,j}) | \xi_{[k-1]}\right]\right\}\\
&{=}&E\left[\nabla\psi_{\mu}(x_k)\right].
\end{eqnarray*}
Thus Lemma \ref{lemma3.3} (a) and (b) hold. Lemma \ref{lemma3.3} (c) can be deduced from the proof of (4.12) in Theorem 2 \cite{Ghadimi}, by noticing the definitions of $\delta_{\mu}^i(x_k)$ and $S^{i-1}_{\mu}(x_{k})$ in \eqref{notation2-2} and \eqref{notation2-3}, respectively.
\end{proof}

In our scheme we update recursively three sequences of points $\{x_k\}_{k=0}^{\infty}$, $\{y_k\}_{k=0}^{\infty}$ and $\{z_k\}_{k=0}^{\infty}$ in $X$. To be specific, given an initial point $x_0\in X$, a fixed smoothing parameter $\mu>0$,
a sequence of mini-batch sizes $\{m_k\}_{k=0}^{\infty}$, and sequences of positive
real numbers $\{\gamma_k\}_{k=0}^{\infty}$, $\{\alpha_k\}_{k=0}^{\infty}$ and $\{\tau_k\}_{k=0}^{\infty}$, then for $k\ge 0$
we recursively obtain
\begin{eqnarray}\label{yk}
y_{k}&=&\arg\min\limits_{y\in X}\left\{ \ell^k(y):=\langle \nabla \Psi_{\mu}^{k},y-x_{k}\rangle+\frac{1}{2 \gamma_{k}}\|y-x_{k}\|^2\right\},\label{yk}\\
z_{k}&=&\arg\min\limits_{x\in X} \left\{\pi^k(x):= \frac{L}{\sigma_d}d(x) +\sum\limits_{i=0}^{k}\alpha_{i}\left[\Psi_{\mu}^i+\langle\nabla \Psi_{\mu}^i,x-x_i\rangle\right]\right\},\label{zk}\\
x_{k+1}&=&\tau_kz_k+(1-\tau_k)y_k. \label{xk+1}
\end{eqnarray}
Let us denote
\begin{eqnarray}\label{BC}
B_i := \dfrac{\sqrt{A_{i}}\sigma^2}{m_iL},\quad C_i:= \dfrac{LA_{i}(1-\sqrt{A_{i}})}{2}\|y_i-x_i\|^2,
\end{eqnarray}
and
\begin{eqnarray}\label{Gamma_k}
\Gamma_k(x):=
E\left[\pi^k(x)\right] + \sum\limits_{i=0}^{k} B_i + E\left[\sum\limits_{i=0}^{k}C_i\right],\quad\quad
\end{eqnarray}
where the expectation is taken with respect to the history of mini-batch samples $\xi_{[k]}$.
By selecting proper positive sequences $\{\gamma_k\}_{k=0}^{\infty}$, $\{\alpha_k\}_{k=0}^{\infty}$ and $\{\tau_k\}_{k=0}^{\infty}$,  we then have the following proposition.

\begin{proposition}\label{Proposition}
Let some sequence $\{\alpha_k\}_{k=0}^{\infty}$ satisfy the condition:
\begin{eqnarray}\label{eq3.15}
\alpha_0\in (0,1], ~~\alpha_k>0,
~~\alpha_{k+1}^2\leq\sqrt{ A_{k+1}},
\end{eqnarray}
where $A_{k}=\sum\limits_{i=0}^{k}\alpha_i$.
Let us choose $\tau_k=\frac{\alpha_{k+1}}{A_{k+1}}$, $\gamma_{k} = \frac{1}{L\sqrt{A_k}}$.
%
Then for any $k\geq0$,
\begin{flalign*}
\quad\quad\quad\quad\quad\quad\quad\quad A_{k}E[\psi_\mu(y_k)]\leq\min \left\{\Gamma_k(x)\ :\ x\in X \right\},\quad\quad\quad\quad\quad\quad(\mathcal{R}_k)
\end{flalign*}
where $\Gamma_k$ is defined in \eqref{Gamma_k} and the expectation is taken with respect to $\xi_{[k]}$.
\end{proposition}
\begin{proof}
For any $k\ge 0$, using Lemma \ref{Lemma2} and the definition of  $\delta_{\mu}(x_k)$ in \eqref{notation2-1},
\begin{eqnarray*}
&&\left\langle\delta_{\mu}(x_k),P_X(x_k,\nabla\Psi_{\mu}^k,\gamma_{k})-P_X(x_k,\nabla\psi_{\mu}(x_k),\gamma_{k})\right\rangle\\
&&~~~\leq ~ \|\delta_{\mu}(x_k)\|\|P_X(x_k,\nabla\Psi_{\mu}^k,\gamma_{k})-P_X(x_k,\nabla\psi_{\mu}(x_k),\gamma_{k})\|\\
&&~\overset{\eqref{Pleqg}}{\leq} \|\delta_{\mu}(x_k)\|\|\nabla\Psi_{\mu}^k-\nabla\psi_{\mu}(x_k)\| =~\|\delta_{\mu}(x_k)\|^2,
\end{eqnarray*}
and
\begin{eqnarray*}
&&E\left[\left\langle\delta_{\mu}(x_k),P_X(x_k,\nabla\psi_{\mu}(x_k),\gamma_{k})\right\rangle\right]\\
&&~=E\left\{E\left[\left\langle\delta_{\mu}(x_k),P_X(x_k,\nabla\psi_{\mu}(x_k),\gamma_{k})\right\rangle |\xi_{[k-1]}\right]\right\}\\
&&~=E\left\{\left\langle E\left[\delta_{\mu}(x_k) |\xi_{[k-1]}\right],P_X(x_k,\nabla\psi_{\mu}(x_k),\gamma_{k})\right \rangle\right\}\\
&&~=E\left\{\left\langle E\left[\left(\nabla\Psi_{\mu}^k-\nabla\psi_{\mu}(x_k)\right) |\xi_{[k-1]}\right], P_X(x_k,\nabla\psi_{\mu}(x_k),\gamma_{k})\right\rangle\right\} = 0.
\end{eqnarray*}
Therefore, for any $k\ge 0$,
\begin{eqnarray}\label{theorem-1-result-1}
&~&A_{k}E\left[\langle\delta_{\mu}(x_k),x_k-y_k\rangle\right]\nonumber\\
&~&~~=\frac{\sqrt{A_{k}}}{L}E\left[\langle\delta_{\mu}(x_k),P_X(x_k,\nabla\Psi_{\mu}^k,\gamma_{k})\rangle\right]\nonumber\\
&~&~~=\frac{\sqrt{A_{k}}}{L}E\left[\langle \delta_{\mu}(x_k),P_X(x_k,\nabla\Psi_{\mu}^k,\gamma_{k})-P_X(x_k,\nabla\psi_{\mu}(x_k),\gamma_{k})\rangle\right]\nonumber\\
&~&~~~~~~~~~+\frac{\sqrt{A_{k}}}{L}E\left[\langle\delta_{\mu}(x_k),P_X(x_k,\nabla\psi_{\mu}(x_k),\gamma_{k})\rangle\right]\nonumber\\
&~&~~ \leq \frac{\sqrt{A_{k}}}{L}E\left[\|\delta_{\mu}(x_k)\|^2\right] \leq B_k,
\end{eqnarray}
where the last inequality comes from Lemma \ref{lemma3.3} (c) and the definition of $B_k$ in (\ref{BC}).

Now we prove the theorem by mathematical induction. For $k=0$, in terms of Lemma \ref{lemma3.3} for the first equality below, we have for any $x\in X$ and $\mu>0$,
\begin{eqnarray*}
&~&A_0E\left[\psi_{\mu}(y_0)\right]\quad\quad\quad\quad\quad\quad\quad\\
&~&~\overset{\eqref{eq3.10}}{\leq} A_0E\left[\psi_{\mu}(x_0)+\langle\nabla\psi_{\mu}(x_0),y_0-x_0\rangle
+\frac{L}{2}\|y_0-x_0\|^2\right]\\
&~&~\overset{\eqref{notation2-1}}{=}A_0E\left[\Psi_{\mu}^0+\langle\nabla\Psi_{\mu}^0,y_0-x_0\rangle+\frac{L}{2}\|y_0-x_0\|^2\right]
+A_0E\left[\langle\delta_{\mu}(x_0),x_0-y_0\rangle\right]\\
&~&~\overset{\eqref{theorem-1-result-1}}{\leq}A_0E\left[\Psi_{\mu}^0+ \ell^0(y_0)\right]
+ C_0 + B_0,
\end{eqnarray*}
where the definition of $\ell^0(y_0)$ can be found in (\ref{yk}).
Using the formula to obtain $y_0$ in (\ref{yk}), we have for any $x\in X$,
\begin{eqnarray}\label{rel1}
\ell^0(y_0) \le \ell^0(x) = \langle \nabla \Psi_{\mu}^0,x-x_0 \rangle + \frac{1}{2\gamma_0} \|x-x_0\|^2.
\end{eqnarray}
Noting that $A_0 = \alpha_0 \in (0,1]$, $\gamma_0 = \frac{1}{L \sqrt{A_0}}$,
and \eqref{eq3.11}, we find the second term in the right-hand side of the equation in (\ref{rel1}) satisfies
\begin{eqnarray}\label{rel2}
\frac{1}{2 \gamma_0} \|x-x_0\|^2 = \frac{L \sqrt{A_0}}{2} \|x-x_0\|^2 \le \frac{L}{\sigma_d} d(x).
\end{eqnarray}
Therefore, by using (\ref{rel1}), (\ref{rel2}), and the definition of $\pi^0(x)$ obtained by (\ref{zk}), we get
\begin{eqnarray*}
A_0\left[\Psi_{\mu}^0 + \ell^0(y_0)\right] \le \frac{L}{\sigma_d}d(x) + \alpha_0 \left[\Psi_{\mu}^0 + \langle \nabla \Psi_{\mu}^0 ,x-x_0\rangle\right] = \pi^0(x).
\end{eqnarray*}
With these observations, we have $(\mathcal{R}_0)$ holds.

Assume $(\mathcal{R}_k)$ holds. By the formula for obtaining $z_k$ in (\ref{zk}),
the optimality condition of the minimization problem in (\ref{zk}) at $z_k$ gives
\begin{eqnarray}\label{argu1}
-\frac{L}{\sigma_d}\langle\nabla d(z_k), x-z_k\rangle\leq
\sum\limits_{i=0}^{k}\alpha_i\left\langle\nabla\Psi_{\mu}^i,x-z_k\right\rangle.
\end{eqnarray}
This, together with the strong convexity of the function $d$ with the convexity parameter $\sigma_d$, yields
\begin{eqnarray}\label{proof-statement-2}
\frac{L}{\sigma_d} \left[d(z_k) + \frac{\sigma_d}{2}\|x-z_k\|^2\right] &\le& \frac{L}{\sigma_d} \left[d(x) - \left\langle \nabla d(z_k),x-z_k \right\rangle\right] \nonumber\\
&\overset{\eqref{argu1}}{\le}& \frac{L}{\sigma_d} d(x) + \sum_{i=0}^k \alpha_i \left\langle \nabla \Psi_{\mu}^i, x-z_k \right\rangle.
\end{eqnarray}
By noting the objective function $\pi_k$ of the minimization problem to find $z_k$ in (\ref{zk}), we know
\begin{eqnarray}\label{pik+1x}
\pi^{k+1}(x) = \pi^{k}(x) + \alpha_{k+1} \left[\Psi_{\mu}^{k+1} + \langle \nabla \Psi_{\mu}^{k+1}, x - x_{k+1} \rangle\right],
\end{eqnarray}
and
\begin{eqnarray}\label{pikx}
\pi^k (x) &=& \frac{L}{\sigma_d} d(x) + \sum_{i=0}^k \alpha_i \left[ \Psi_{\mu}^i + \langle \nabla \Psi_{\mu}^i, x-x_i\rangle\right] \nonumber\\
           &=&  \frac{L}{\sigma_d} d(x) +  \sum_{i=0}^k \alpha_i \langle \nabla \Psi_{\mu}^i, x-z_k \rangle +
                   \sum_{i=0}^k \alpha_i \left[\Psi_{\mu}^i + \langle \nabla \Psi_{\mu}^i, z_k-x_i\rangle \right] \nonumber\\
            &\overset{\eqref{proof-statement-2}}{\ge}& \frac{L}{\sigma_d}\left[d(z_k) + \frac{\sigma_d}{2}\|x-z_k\|^2\right] + \sum_{i=0}^k \alpha_i \left[\Psi_{\mu}^i + \langle \nabla \Psi_{\mu}^i, z_k-x_i\rangle \right] \nonumber\\
            &=& \pi^k(z_k) + \frac{L}{2}\|x-z^k\|^2.
\end{eqnarray}

Let us denote
\begin{eqnarray}
T^{k+1} &:=& \alpha_{k+1} \left[\Psi_{\mu}^{k+1} + \langle \nabla \Psi_{\mu}^{k+1}, x - x_{k+1} \rangle\right], \label{Tk+1}\\
t^{k+1} &:=& \alpha_{k+1}  \left[\psi_\mu(x_{k+1})+\langle \nabla
        \psi_\mu(x_{k+1}), x-x_{k+1}\rangle\right]. \label{tk+1}
\end{eqnarray}
It is clear that $E\left[T^{k+1}\right] = E\left[t^{k+1}\right]$ by Lemma \ref{lemma3.3}. Combining \eqref{pik+1x} and \eqref{pikx}, we have
\begin{eqnarray}\label{proof-statement-3}
\pi^{k+1}(x) \ge \pi^k(z_k) + \frac{L}{2} \|x-z_k\|^2 + T^{k+1}.
\end{eqnarray}
Therefore, for any $x\in X$,
\begin{eqnarray}\label{Gammak+1x}
\Gamma_{k+1}(x)
 &=& E\left[\pi^{k+1}(x)\right] + \sum_{i=0}^{k+1} B_i + E\left[\sum_{i=0}^{k+1} C_i\right]  \nonumber\\
                &\overset{\eqref{proof-statement-3}}{\ge}& E\left[\pi^k(z_k)+\frac{L}{2}\|x-z^k\|^2+ T^{k+1}\right] + \sum_{i=0}^{k+1} B_i + E\left[\sum_{i=0}^{k+1} C_i\right] \nonumber\\
                &=& \Gamma_k(z_k)  + E\left[T^{k+1}\right] + \frac{L}{2} E\left[\|x- z_k\|^2\right] + B_{k+1} + E\left[C_{k+1}\right].~~
\end{eqnarray}
For the summation of the first two terms in (\ref{Gammak+1x}), we have
\begin{eqnarray}\label{succeed-1}
& & \Gamma_k(z_k)+ E\left[T^{k+1}\right] \nonumber\\
 & &  \quad \geq    A_k E\left[\psi_\mu(y_k)\right] + E\left[t^{k+1}\right] \nonumber\\
 & & \quad  \geq   A_k E\left[\psi_\mu(x_{k+1}) + \langle \nabla \psi_{\mu}(x_{k+1}), y_k - x_{k+1}\rangle \right] + E\left[t^{k+1}\right] \nonumber\\
& & \quad = A_{k+1}E\left[\psi_\mu(x_{k+1})\right] + E\left[\langle\nabla\psi_\mu(x_{k+1}),A_{k}y_k-A_{k+1}x_{k+1}+\alpha_{k+1}x\rangle\right].\quad\quad\quad
\end{eqnarray}
Here the first inequality comes form the assumption $({\cal R}_k)$ holds, and the fact that $E[T^{k+1}] = E[t^{k+1}]$,
and the second inequality is obtained from the convexity of $\psi_{\mu}$, and the first equality is obtained by
using the definition of $t^{k+1}$ in (\ref{tk+1}), and the choice of $A_k$ such that $A_{k+1} = A_k + \alpha_{k+1}$.
Furthermore, for the second term of $\langle \cdot , \cdot \rangle$ in (\ref{succeed-1}), by noting
$\tau_k=\frac{\alpha_{k+1}}{A_{k+1}}$, $1-\tau_k = \frac{A_k}{A_{k+1}}$, and $x_{k+1}=\tau_kz_k+(1-\tau_k)y_k$ in (\ref{xk+1}),
we have
\begin{eqnarray}\label{succeed-2}
A_{k}y_k-A_{k+1}x_{k+1}+\alpha_{k+1}x
&=&A_{k}y_k-\alpha_{k+1}z_{k}-A_{k+1}\frac{A_{k}}{A_{k+1}}y_k+\alpha_{k+1}x \nonumber\\
&=&\alpha_{k+1}(x-z_{k}) = A_{k+1}\tau_k (x-z_k).
\end{eqnarray}
In view of $\tau_k = \frac{\alpha_{k+1}}{A_{k+1}}$ and $\alpha_{k+1}^2\leq\sqrt{ A_{k+1}}$, we know
\begin{eqnarray}\label{coefficient}
\frac{1}{ A_{k+1}}\geq\tau_k^2\sqrt{ A_{k+1}}.
\end{eqnarray}
Using \eqref{Gammak+1x} -- \eqref{coefficient}, we find
\begin{eqnarray}\label{succeed-1'}
&&\Gamma_{k+1}(x)\nonumber\\
&&~~\ge A_{k+1}E\left[\psi_{\mu}(x_{k+1})\right] + A_{k+1} \tau_k E\left[\langle \nabla \psi_{\mu}(x_{k+1}),x-z_k \rangle\right] + \frac{L}{2} E\left[\|x-z_k\|^2\right] \nonumber\\
&&~~~~~~~~+B_{k+1} + E\left[C_{k+1}\right] \nonumber \\
&&~\overset{\eqref{coefficient}}{\geq} A_{k+1}E\left[\psi_\mu(x_{k+1}) +\tau_k\langle\nabla\psi_\mu(x_{k+1}),x-z_k\rangle+\frac{L\tau_k^2\sqrt{A_{k+1}}}{2}\|x-z_k\|^2\right] \nonumber\\
&&~~~~~~~~+B_{k+1}+E\left[C_{k+1}\right].
\end{eqnarray}

For any $x\in X$, let $y(x,k) := \tau_k x + (1- \tau_k)y_k $. By using $x_{k+1}= \tau_k z_k + (1- \tau_k) y_k $ in (\ref{xk+1}),
we get
\begin{eqnarray*}
y(x,k) - x_{k+1} = \tau_k (x- z_k),
\end{eqnarray*}
and consequently
\begin{eqnarray}\label{succeed-2'}
 \ell^{k+1}\left(y(x,k)\right) &=& \langle\nabla\Psi_\mu^{k+1},y(x,k)-x_{k+1}\rangle
+\frac{1}{2\gamma_{k+1}}\|y(x,k)-x_{k+1}\|^2 \nonumber \\
&=& \tau_k\langle\nabla\Psi_\mu^{k+1},x-z_k\rangle
+\frac{L\tau_k^2\sqrt{A_{k+1}}}{2}\|x-z_k\|^2.
\end{eqnarray}
Note that $x$ and $z_k$ are deterministic if the history of mini-batch samples $\xi_{[k]}$ is given.
Thus
\begin{eqnarray}
E[\langle \nabla \Psi_{\mu}^{k+1},x-z_k \rangle] &=& E \{ E[\langle \nabla \Psi_{\mu}^{k+1},x-z_k  \rangle \mid_{\xi_{[k]}}] \} \nonumber\\
   &=& E\{ \langle E[ \nabla \Psi_{\mu}^{k+1} ] \mid_{\xi_{[k]}}, x-z_k \rangle]\} \nonumber\\
   &=& E[ \langle \nabla \psi_{\mu}(x_{k+1}),x-z_k \rangle ]. \label{succeed-3'}
\end{eqnarray}
Therefore, by combining \eqref{succeed-1'}, \eqref{succeed-2'}, \eqref{succeed-3'}
and using $E[\psi_{\mu}(x_{k+1})] = E[\Psi_{\mu}^{k+1}]$ according to Lemma \ref{lemma3.3} (a), we find
\begin{eqnarray}\label{eq3.8}
\Gamma_{k+1}(x)\geq A_{k+1}E\left[\Psi_\mu^{k+1} +\ell^{k+1}(y(x,k))\right]
+B_{k+1}+E\left[ C_{k+1}\right].
\end{eqnarray}
On the other hand,
\begin{eqnarray}
& &A_{k+1}E\left[\psi_\mu(y_{k+1})\right]\nonumber\\
& &~\overset{\eqref{eq3.10}}{\leq} A_{k+1}E\left[\psi_\mu(x_{k+1})+\langle\nabla\psi_\mu(x_{k+1}),y_{k+1}-x_{k+1}\rangle
+\frac{L}{2}\|y_{k+1}-x_{k+1}\|^2\right]\nonumber\\
& &~~~= A_{k+1}E\left[\Psi_\mu^{k+1}+\langle\nabla\Psi_\mu^{k+1} -\delta_\mu(x_{k+1}),y_{k+1}-x_{k+1}\rangle
+\frac{L}{2}\|y_{k+1}-x_{k+1}\|^2\right]\nonumber\\
& &~~~= A_{k+1}E \left[\Psi_{\mu}^{k+1} + \ell^{k+1}(y_{k+1}) \right] + A_{k+1} E\left[\langle \delta_{\mu}(x_{k+1}), x_{k+1} - y_{k+1}\rangle\right] \nonumber\\
& &~~~~~~~\quad\quad + E\left[ A_{k+1}(\frac{L}{2} - \frac{1}{2 \gamma_{k+1}}) \|y_{k+1}-x_{k+1}\|^2\right] \nonumber\\
& &~\overset{\eqref{theorem-1-result-1}}{\leq} A_{k+1}E\left[\Psi_\mu^{k+1}+\ell^{k+1}(y(x,k))\right]
+B_{k+1}+E\left[ C_{k+1}\right],\label{eq3.7}
\end{eqnarray}
where the first equality is deduced by using Lemma \ref{lemma3.3} (a), and the definition of $\delta_{\mu}(x_{k+1})$, and the last inequality is obtained due to $y_{k+1}=\arg\min\limits_{y\in X}\left\{ \ell^{k+1}(y)\right\}$, \eqref{theorem-1-result-1}, $\gamma_{k+1} = \frac{1}{L\sqrt{A_{k+1}}}$, and the definition for $C_{k+1}$ in (\ref{BC}).
It is easy to see from \eqref{eq3.8} and \eqref{eq3.7} that $A_{k+1}E[\psi_\mu(y_{k+1})]\leq \min\{\Gamma_{k+1}(x)\ :\ x\in X\}$. That is, (${\cal R}_{k+1}$) holds.

Therefore, by using mathematical induction, the relation $(\mathcal{R}_k)$ holds for any $k\geq0$.
\end{proof}

Clearly, there are many ways to satisfy the conditions \eqref{eq3.15}. Here we choose the special $\alpha_k$ and batch sizes
$m_k$ for $k\ge 0$ as follows in order to guarantee the optimal iteration complexity of the min-batch stochastic Nesterov's smoothing method we propose later.

\begin{corollary}\label{corollary}
For $k\geq 0$ define $\alpha_{k}\equiv\frac{1}{2}$. Then
\begin{eqnarray*}
A_{k}=\sum\limits_{i=0}^{k}\alpha_{i}=\frac{k+1}{2},\ \tau_k=\frac{\alpha_{k+1}}{A_{k+1}}=\frac{1}{k+2},
\end{eqnarray*}
and the conditions in \eqref{eq3.15} are satisfied. If in addition,  the batch sizes $m_k \equiv m\geq 1$ for $k = 0,\ldots, N$, then
\begin{eqnarray}\label{cor-1}
E[\psi_\mu(y_N)]&\leq&\min\limits_{x\in X}\left\{\sum\limits_{i=0}^{N}\frac{1}{N+1}E\left[\psi_\mu(x_i)+\langle\nabla \psi_\mu(x_i),x-x_i\rangle\right]\right\}\nonumber\\
&~&~\quad +\ \frac{(6-\sqrt{2})LD}{2(N+1)\sigma_{d}}+\frac{\sqrt{2(N+1)}\sigma^2}{mL}.
\end{eqnarray}
\end{corollary}
\begin{proof}
Indeed, $\alpha_{k+1}^2 = \frac{1}{4} \le \sqrt{\frac{k+2}{2}}= \sqrt{A_{k+1}}$ for all $k\ge 0$, and consequently conditions in
 \eqref{eq3.15} are satisfied. For the special choices of $\alpha_k$ and $m_k$ of this corollary, we immediately get $A_i = \frac{i+1}{2}$,
\begin{eqnarray*}
\sum_{i=0}^N B_i = \sum_{i=0}^N \frac{\sqrt{A_i }\sigma^2}{m L}
                 \le \frac{(N+1)\sqrt{N+1} \sigma^2}{\sqrt{2}mL},
\end{eqnarray*}
and
\begin{eqnarray*}
\sum_{i=1}^N C_i &=& \sum_{i=0}^N \frac{L A_i(1-\sqrt{A_i})}{2}\|y_i-x_i\|^2 \\
                 &\le&   \frac{L A_0 (1-\sqrt{A_0})}{2} \|y_0 -x_0\|^2
                 \le  \frac{(2-\sqrt{2})L D}{4\sigma_d},
\end{eqnarray*}
where the first inequality holds, since $1-\sqrt{A_i} \le 0$ for all $i\ge 1$, and the second inequality holds, since
\eqref{eq3.11} and \eqref{eq3.12} implies
\begin{eqnarray*}
\frac{1}{2} \sigma_d \|y_0 - x_0\|^2 \le d(y_0) \le D.
\end{eqnarray*}
 According to Proposition \ref{Proposition}, $({\cal R}_N)$ holds.
Moreover, by using \eqref{eq3.12}, we have
\begin{eqnarray}{\label{in1}}
\frac{2}{N+1}\frac{L}{\sigma_d} d(x) \le 4 \frac{LD}{2\sigma_d (N+1)}.
\end{eqnarray}
Because $x_i$ relates only to $\xi_{[i-1]}$, and $\nabla \Psi_{\mu}^i$ relates to $\xi_{[i]}$,
we get by Lemma \ref{lemma3.3} that
\begin{eqnarray}\label{e2}
E[\Psi_{\mu}^i] = E[\psi_{\mu}(x_i)],\quad \mbox{and}\quad E[\langle \nabla \Psi_{\mu}^i, x-x_i \rangle ] = E[\langle \nabla \psi_{\mu}(x_i), x-x_i \rangle ],
\end{eqnarray}
 Therefore,
\begin{eqnarray*}
E[\psi_{\mu}(y_N)] &~\overset{({\cal R}_N)}{\le}& \min_{x\in X}\frac{2}{N+1} E \left\{ \frac{L}{\sigma_d} d(x) +
\sum_{i=0}^N \frac{1}{2}[\Psi_{\mu}^i + \langle \nabla \Psi_{\mu}^i,x-x_i \rangle]  \right\} \\
 & & \quad\quad + \frac{\sqrt{2(N+1)}\sigma^2}{mL} + \frac{(2-\sqrt{2})LD}{2 \sigma_d (N+1)},
\end{eqnarray*}
and consequently  we get
\eqref{cor-1} as we desired,
by further noting (\ref{in1}), and the equalities in (\ref{e2}).
\end{proof}

We are ready to give the mini-batch stochastic Nesterov's smoothing method as follows.
%

\vspace{3mm}\hspace{-5mm}\textbf{Algorithm 1}~Given initial point $x_0\in X$, iteration limit $N$, the batch sizes $m_{k}\equiv m >0$ for all $k$.\vspace{2mm}\\
For $k=0:{N}$ do\vspace{1mm}\\
\vspace{1mm}
\hspace{5mm}1. Call the $\cal SO$ $m_{k}$ times to obtain $\Psi_{\mu}(x_k, \xi_{k,i})$ and $\nabla \Psi_{\mu}(x_k, \xi_{k,i})$, $i=1,\cdots$, \\
\vspace{1mm}
\hspace{9mm}$m_{k}$, and set\vspace{1mm}\\
\vspace{1mm}
\hspace{17mm}$\Psi_{\mu}^k=\frac{1}{m_{k}}\sum\limits_{i=1}^{m_{k}} \Psi_{\mu}(x_k, \xi_{k,i})$, $\nabla \Psi_{\mu}^k=\frac{1}{m_{k}}\sum\limits_{i=1}^{m_{k}}\nabla \Psi_{\mu}(x_k, \xi_{k,i})$.\vspace{1mm}\\
\vspace{1mm}
\hspace{5mm}2. Find $y_k=\arg\min\limits_{y\in X}\left\{\langle \nabla \Psi_{\mu}^k,y-x_k\rangle+\frac{L\sqrt{k+1}}{2\sqrt{2}}\|y-x_k\|^2\right\}$.\vspace{1mm}\\
\vspace{1mm}
\hspace{5mm}3. Find $z_{k}=\arg\min\limits_{x\in X} \left\{\frac{L}{\sigma_d}d(x) +\sum\limits_{i=0}^{k}\frac{1}{2}\left[\Psi_{\mu}^i+\langle\nabla \Psi_{\mu}^i,x-x_i\rangle\right]\right\}$.\vspace{1mm}\\
\vspace{1mm}
\hspace{5mm}4. Set $x_{k+1}=\frac{1}{k+2}z_k+\frac{k+1}{k+2}y_k$.\\
\textbf{Output:}~$y_N$.

\vskip 3mm
Recall that $u_{\mu}(x)$ is the unique optimal solution of
\begin{eqnarray*}
\max_{u\in U} \left\{\langle Ax, u \rangle - Q(u) - \mu \omega(u)   \right\}.
\end{eqnarray*}
We will use $u_{\mu}(x)$ in the convergence theorem of the mini-batch stochastic Nesterov's smoothing method.

\begin{theorem}\label{convergence}
Let us apply Algorithm 1 to the smooth problem $\min\limits_{x\in X}\psi_\mu(x)$ with the following value of smoothing parameter:
\begin{eqnarray}\label{eq3.23}
\mu=\mu(N)=\frac{\|A\|^2\sqrt{(6-\sqrt{2})mD}}{\sqrt{2(N+1)\sigma_{d} \sigma_\omega}
	\sqrt{m\|A\|^2\Omega+\sqrt{2(N+1)}\sigma_\omega\sigma^2}},
\end{eqnarray}
where the batch size $m$ is given by
\begin{eqnarray}\label{eq3.34}
	m= \left\lceil \frac{\sqrt{2}\sigma^2\sigma_\omega\sqrt{N+1}}{\|A\|^2\Omega} \right\rceil.
\end{eqnarray}
Here $\lceil \beta \rceil$ means the smallest integer that is no less than a given $\beta>0$.
Then after $N$ iterations we can generate the approximate solution $\widehat x = y_N$ to the original problem \eqref{orip}
that satisfy the following inequality:
\begin{eqnarray}\label{eq3.24}
0\leq E\left[\psi(\widehat{x})-\phi(\widehat{u})\right]
&\leq& \frac{2 \|A\|\sqrt{(6-\sqrt{2})D\Omega}}{\sqrt{(N+1)\sigma_{d} \sigma_\omega}}
+\frac{(6-\sqrt{2}) L_fD}{(N+1)\sigma_{d}},
\end{eqnarray}
where
\begin{eqnarray}\label{uN}
\widehat{u}=\sum\limits_{i=0}^{N}\frac{1}{N+1}u_{\mu}(x_i)\in U.
\end{eqnarray}
Consequently, the iteration complexity of finding an $\epsilon$-approximate solution to the original problem \eqref{orip} does not exceed
\begin{eqnarray}\label{complexity}
N+1=\left\lceil\frac{4(6-\sqrt{2})D\Omega\|A\|^2}{\sigma_{d} \sigma_\omega}\cdot\frac{1}{\epsilon^2}
+\frac{2(6-\sqrt{2})L_fD}{\sigma_{d}}\frac{1}{\epsilon}\right\rceil.
\end{eqnarray}
\end{theorem}
\begin{proof}
Let us fix an arbitrary $\mu > 0$. In view of Corollary \ref{corollary}, we find
\begin{eqnarray}
E\left[\psi_\mu(\widehat{x})\right]&\leq&\min\limits_{x\in X}\left\{\sum\limits_{i=0}^{N}\frac{1}{N+1}E\left[\psi_\mu(x_i)+\langle\nabla \psi_\mu(x_i),x-x_i\rangle\right]\right\}\nonumber\\
&~&~\quad \quad +\frac{(6-\sqrt{2})LD}{2(N+1)\sigma_{d}}+\frac{\sqrt{2(N+1)}\sigma^2}{mL}.\label{eq3.25}
\end{eqnarray}
Recall that $u_{\mu}(x)$ is the unique optimal solution of \eqref{eq2.6}. We have for any $x\in X$,
\begin{eqnarray*}
 h_{\mu}(x)=\langle Ax,u_{\mu}(x)\rangle-Q\left(u_{\mu}(x)\right)-\mu\omega(u_{\mu}(x)),
\end{eqnarray*}
and by Lemma \ref{lemm2.1}
\begin{eqnarray*}
\langle\nabla h_{\mu}(x),x\rangle = \langle A^Tu_{\mu}(x),x\rangle,
\end{eqnarray*}
which, together with (\ref{eq2.8}), imply
\begin{eqnarray}\label{equa1}
h_{\mu}(x) - \langle \nabla h_{\mu}(x),x \rangle = -Q\left(u_{\mu}(x)\right)-\mu\omega(u_{\mu}(x)) \le -Q\left(u_{\mu}(x)\right).
\end{eqnarray}
By (\ref{uN}) and the convexity of $Q(u)$, we know that
\begin{eqnarray}\label{Qc}
Q(\widehat{u}) = Q\left(\sum\limits_{i=0}^{N} \frac{1}{N+1} u_{\mu}(x_i)\right) \le \sum\limits_{i=0}^{N}\frac{1}{N+1}Q(u_{\mu}(x_i)).
\end{eqnarray}
Let us denote
\begin{eqnarray}\label{tildex}
\widetilde{x}=\arg\min \left\{f(x)+\langle A^T\widehat{u},x\rangle \ : x\in X\right\}.
\end{eqnarray}
We then have the first term in the right-hand side of (\ref{eq3.25}) satisfies
\begin{eqnarray*}
&~&\min\limits_{x\in X}\left\{\sum\limits_{i=0}^{N}\frac{1}{N+1}E\left[\psi_\mu(x_i)+\langle\nabla \psi_\mu(x_i),\widetilde x-x_i\rangle\right]\right\}  \nonumber\\
&~&~\quad \le \sum\limits_{i=0}^{N}\frac{1}{N+1}E\left[\psi_\mu(x_i)+\langle\nabla \psi_\mu(x_i),\widetilde  x-x_i\rangle\right]\\
&~&~\quad \overset{\eqref{smooth}}{=}\sum\limits_{i=0}^{N}\frac{1}{N+1}E[f(x_i)+\langle\nabla f(x_i),\widetilde x-x_i\rangle + \langle \nabla h_{\mu}(x_i),\widetilde  x\rangle   \\
& &\quad\quad \quad\quad\quad\quad + h_{\mu}(x_i)-\langle\nabla h_{\mu}(x_i),x_i\rangle]\\
&~&~\quad \overset{\eqref{equa1}}{\leq} \sum\limits_{i=0}^{N}\frac{1}{N+1}E\left[ f(\widetilde  x)+\langle
      A^Tu_{\mu}(x_i),\widetilde  x\rangle-Q\left(u_{\mu}(x_i)\right)\right]\\
&~&~\quad \overset{\eqref{Qc}}{\leq}  E\left[f(\widetilde{x})+\langle
       A^T\widehat{u},\widetilde{x}\rangle-Q(\widehat{u})\right]\\
&~&~\quad \overset{\eqref{tildex}}{=} E\left[-Q(\widehat{u})+\min\limits_{x\in X}\left\{f(x)+\langle A^T\widehat{u},x\rangle\right\}\right]\\ &~&~\quad \overset{\eqref{eq2.5}}{=} E\left[\phi(\widehat{u})\right],
\end{eqnarray*}
where the first inequality also employs the convexity of $f$ and Lemma \ref{lemm2.1}.
Hence, the above inequality, together with \eqref{eq2.7} and \eqref{eq3.25} yields
\begin{eqnarray*}
E\left[\psi(\widehat{x})-\phi(\widehat{u})\right]-\mu\Omega\overset{\eqref{eq2.7}}{\leq} E\left[\psi_\mu(\widehat{x})-\phi(\widehat{u})\right]\leq \frac{\sqrt{2(N+1)}\sigma^2}{mL}+\frac{\left(6-\sqrt{2}\right)LD}{2(N+1)\sigma_{d}}.
\end{eqnarray*}
By Lemma \ref{lemm2.1}, the gradient of $\psi_\mu(x)$ is Lipschitz continuous with the constant
\begin{eqnarray*}
L = L_f+ L_{h_{\mu}} = L_f+\frac{1}{\mu \sigma_\omega}\|A\|^2.
\end{eqnarray*}
Then
\begin{eqnarray}
0&\leq& E\left[\psi(\widehat{x})-\phi(\widehat{u})\right] \nonumber\\
&\leq&\mu\Omega+\frac{\sqrt{2(N+1)}\sigma^2}{mL}+\frac{(6-\sqrt{2})LD}{2(N+1)\sigma_{d}} \nonumber \\
&=&\mu\Omega+\frac{\sqrt{2(N+1)}\sigma^2}{m\left(L_f+\frac{1}{\mu \sigma_\omega}\|A\|^2 \right)}+\frac{(6-\sqrt{2})\left(L_f+\frac{1}{\mu \sigma_\omega}\|A\|^2 \right)D}{2(N+1)\sigma_{d}} \nonumber\\
&\le&\mu\left[\Omega+\frac{\sqrt{2(N+1)}\sigma_\omega\sigma^2}{m\|A\|^2}\right]
+\frac{1}{\mu}\frac{(6-\sqrt{2})\|A\|^2D}{2(N+1)\sigma_{d} \sigma_\omega}
+\frac{(6-\sqrt{2})L_fD}{2(N+1)\sigma_{d}}. \label{havemum}
\end{eqnarray}
Note that $a+b\geq 2\sqrt{ab}$ for any positive real numbers $a,b$, and the equality holds
 if and only if $a=b$. Thus by choosing $\mu = \mu(N)$ as in (\ref{eq3.23}), we find the first two terms of (\ref{havemum})
    are equal, and the right-hand side of this inequality in $\mu$ is minimized.
Letting the batch size $m$ be given by \eqref{eq3.34},
 we then get \eqref{eq3.24} from \eqref{havemum}.

Letting the right hand side of \eqref{eq3.24} smaller than $\epsilon$,
we then easily know that the iteration complexity of finding an $\epsilon$-approximate solution to the original problem \eqref{orip}
does not exceed \eqref{complexity}.
\end{proof}

\begin{remark}
%
By \eqref{complexity}, we can conclude that Algorithm 1 has the optimal iteration complexity $O(\frac{1}{\epsilon^2})$ to get an $\epsilon$-approximate solution. The special choices of $\alpha_k \equiv \frac{1}{2}$, $\gamma_k = \frac{\sqrt{2}}{L \sqrt{k+1}}$, and $m_k \equiv m$ in Algorithm 1 are important to guarantee the optimal iteration complexity of the mini-batch stochastic Nesterov's smoothing method.
\end{remark}

\section{Application in support vector machine}
Support vector machine (SVM) is a popular machine learning method for classification \cite{Noble,Huang,Rodriguez}. After solving certain optimization problems on known data samples with labels by some method, the parameters of the classifier are determined. The determined classifier is then used to predict the labels of new data samples without label information.

In order to evaluate the performance of a certain SVM model associated with a certain algorithm, users often divide the known data samples with label information into two parts: training data and testing data. Training data means a set of data samples used for learning, which is to fit the parameters of the classifier. Testing data refers to a set of data samples used only to assess the performance of the classifier. After solving optimization problems on training data by some method, users can apply decision functions to predict the labels of testing data. Let $z_1, \ldots, z_J$ be the testing data and $\bar {y}_1, \ldots, \bar{y}_J$ be the predicted labels. If the true labels of testing data are known and denoted as $y_1, \ldots, y_J$, $\mathcal{S}=\{i\in\{1,\ldots,J\}:\bar {y}_i=y_i\}$, we evaluate the prediction results by the following measure:
\begin{eqnarray*}
\mathrm{Accuracy~(Acc)}
&=&\frac{\mathrm{the~number~of~correctly~predicted~data\ (\sharp \mathcal{S})}}{\mathrm{the~number~of~total~testing~data}\ (J)}\times 100\%.
\end{eqnarray*}

In this section, we consider an application of our mini-batch stochastic Nesterov's smoothing method on a stochastic nonsmooth convex model of SVM for binary classification as follows \cite{Shivaswamy}
\begin{eqnarray}\label{SVM-ori}
\begin{array}{ll}
\min\limits_{x}&\lambda_1x^T\Sigma x+ E\left[\max\left\{0,1-y_{\xi}\langle x,z_{\xi}\rangle\right\}\right] \\
\rm{s.t.}&\|x\|^2\leq t.
\end{array}
\end{eqnarray}
Here, $w_{\xi} = (z_{\xi}^T,y_{\xi})^T\in \mathds{R}^{n+1}$ is the random vector, $E$ is the expectation with respect to $\xi$, $\Sigma$ is the covariance matrix of the random vector $z_{\xi}$,
and $t>0$ is a given parameter. Denote by $NS$ the number of the training samples. Then
\begin{eqnarray*}
\Sigma=\frac{1}{NS}\sum\limits_{i=1}^{NS}z_{\xi}^i {z_{\xi}^i}^T-\frac{1}{NS^2}\sum\limits_{i=1}^{NS}z_{\xi}^i\sum\limits_{j=1}^{NS}{z_{\xi}^j}^T,
\end{eqnarray*}
where $z_{\xi}^i$ denotes the $i$-th training sample. This formulation of $\Sigma$ comes from \cite{Shivaswamy}.

For this model, $E\left[\max\left\{0,1-y_{\xi}\left\langle x,z_{\xi}\right\rangle\right\}\right]$ is the nonsmooth term which itself involves the expectation operator. It is not easy to get the proximal mapping for the max operator of an affine function and zero even if only one random vector $w_{\xi}$ is selected.

Let $X=\left\{x\in \mathds{R}^n \ :\ \|x\|^2\leq t\right\}$. We can reformulate the above model as
\begin{eqnarray}\label{eq5.30}
\min\limits_{x\in X}\ \{ \psi(x) =  f(x)+h(x)\},
\end{eqnarray}
where $f(x)=\lambda_1x^T\Sigma x$, and
\begin{eqnarray}
h(x)=E\left[\max\left\{0,1-y_{\xi}z_{\xi}^Tx\right\}\right]
= E\left[\max\limits_{0\leq u\leq 1}\left\{u(1-y_{\xi}z_{\xi}^T x)\right\}\right].\label{eq5.31}
\end{eqnarray}
The Lipschitz constant for $f$ is $L_f=2\lambda_1\|\Sigma\|$, where $\|M\|=\lambda_{\max}(M)$, and $\lambda_{\max}(M)$ means the maximal eigenvalue of $M$. Let $A_{\xi}=-y_{\xi}z_{\xi}^T$, $A=E[A_{\xi}]$, $Q(u) = -u$ and $U = \{u\in \mathds{R}\ : 0\le u \le 1\}$. Then $h(x)$ can be written in the form of \eqref{eq2.6} as
\begin{eqnarray*}
h(x) = \max_{u\in U}\left\{\langle Ax, u \rangle -Q(u) \right\}.
\end{eqnarray*}

Hence our mini-batch stochastic Nesterov's smoothing method is suitable to solve the SVM model defined in \eqref{SVM-ori}.
In our numerical experiment,
we choose the Euclidean distance as the prox-function, that is,
\begin{eqnarray}\label{eq5.32}
d(x)=\frac{1}{2}\sum\limits_{i=1}^{n}x_i^2,~~~\omega(u)=\frac{1}{2}u^2.
\end{eqnarray}
Thus we get $\sigma_{d}=1$, $\sigma_{\omega}=1$, $D=\max\limits_{x\in X} d(x)
=\frac{t}{2}$, and $\Omega=\max\limits_{u\in U} \omega(u)=\frac{1}{2}$.

In view of \eqref{eq5.31}, the smoothing function of $h(x)$ is
\begin{eqnarray*}
h_{\mu}(x)&=&E\left[\max\limits_{u\in U}\left\{u(1-y_{\xi}z_{\xi}^T x)-\frac{\mu}{2}u^2\right\}\right]\\
&=&\left\{ \begin{aligned}
&0,~~~~~~~~~~~~~~~~~~~~~~~~~~~~~~~~~~~y_{\xi}z_{\xi}^Tx > 1,\\
&E\left[\frac{(1-y_{\xi}z_{\xi}^T x)^2}{2\mu}\right],~~~~~~~~~~~~1-\mu\leq y_{\xi}z_{\xi}^T x\leq1,\\
&E\left[1-y_{\xi}z_{\xi}^Tx-\frac{\mu}{2}\right],~~~~~~~~~~~y_{\xi}z_{\xi}^Tx < 1-\mu.
\end{aligned}\right.
\end{eqnarray*}

Given the tolerance $\epsilon>0$, we determine the number of iteration $N$ needed to get the $\epsilon$-approximate solution of \eqref{SVM-ori} according to \eqref{complexity} of Theorem \ref{convergence}, and consequently find the corresponding smoothing factor $\mu$ by \eqref{eq3.23} of Theorem \ref{convergence}. The Nesterov's smoothing method
in fact uses the iterates on solving the smooth problem
\begin{eqnarray}\label{SVM-smooth}
\min_{x\in X}\  \{ \psi_{\mu}(x)=f(x) + h_{\mu}(x)\},
\end{eqnarray}
with the fixed $\mu$, and outputs $y_N$ as the computed $\epsilon$-approximate solution of the original nonsmooth SVM model in \eqref{SVM-ori}.

We use the $NS$ training samples to estimate the parameters, $\|A\|$, $\sigma^2$ and $L$. Specifically, we get $\|A\|=\lambda_{\max}(A)=\lambda_{\max}(E[A_{\xi}])$, where the expectation is taken with respect to $\xi$ in the training data. Then, an estimation of $L_{h_{\mu}}$ is obtained by \eqref{eq2.9}, and the Lipschitz constant for $\psi_{\mu}(x)$ is $L=L_f+L_{h_{\mu}}$. We follow the way of estimating the parameter $\sigma^2$ as in \cite{Ghadimi}. Using the training samples, we compute the stochastic gradients of the objective function $\lceil NS/100\rceil$ times at 100 randomly selected points and then take the average of the variances of the stochastic gradients for each point as an estimation of $\sigma^2$.

The existing SA type methods can not solve the nonsmooth SVM model defined in \eqref{SVM-ori} with guaranteed convergence, because of the relatively complex nonsmooth term that leads to the difficulty of obtaining its proximal operator.
For comparison, we apply the existing SA methods to solve the smooth counterpart in \eqref{SVM-smooth}, including the randomized stochastic projected gradient (RSPG) method, the two-phase RSPG (2-RSPG) method and its variant 2-RSPG-V method in \cite{Ghadimi}, as well as the mini-batch mirror descent SA (M-MDSA) method. The M-MDSA method is a mini-batch version of the MDSA method with constant stepsize policy in \cite{Nemirovski}. Such modification improves the computational speed significantly compared with the original MDSA in \cite{Nemirovski} as pointed out in \cite{Ghadimi}. The batch sizes of the M-MDSA method are set to be the same as that for our MSNS method.

It's worth noting that either the 2-RSPG method, or the 2-RSPG-V method includes two phases -- the optimization phase and the post-optimization phase. Each method generates several candidate outputs in the optimization phase, and the final output is selected from these candidate outputs according to some rules in the post-optimization phase.

We denote by $N_{total}$ the total number of calls to the $\cal SO$. We set $N_{total}$ to be the same for different methods. Since the RSPG, 2-RSPG and 2-RSPG-V methods are randomized SA methods, they stop randomly before up to the maximum number of the $\cal SO$ calls. We then set $N_{total}$ for the three methods to be the maximum number of calls to the $\cal SO$ (for the 2-RSPG and 2-RSPG-V methods, it refers to the maximum number of calls to the $\cal SO$ in the optimization phase).

We provide the numerical experiments on both synthetic datasets and real datasets. Our experiments were performed in MATLAB R2018b on a laptop with 1 dualcore 2.4 GHz CPU and 8 GB of RAM.

\subsection{Synthetic datasets}

Let the training data be given by $T=\{(z_\xi,y_\xi)\}$. Here, we assume that the feature vector $z_\xi$ is drawn from standard normal distribution with approximately $10\%$ nonzero elements. We randomly generate a vector $\bar{x}\in X$. Using this $\bar{x}$, we determine the label $y_\xi \in \{-1, 1\}$ to be
\begin{eqnarray*}
y_\xi&=&\left\{ \begin{aligned}
1,~~~~~~~~~~&\langle\bar{x},z_\xi\rangle\geq0,\\
-1,~~~~~~~~~&\langle\bar{x},z_\xi\rangle<0.
\end{aligned}\right.
\end{eqnarray*}
The testing data contains $K = 50000$ samples.
The smoothing parameter $\mu$ and the batch size $m$ are set according to \eqref{eq3.23} and \eqref{eq3.34}, respectively. The other parameters are set to be $\lambda_1 = 0.5$, $t=10$. The dimensions of the problems are set to be $n = 500$ and $1000$, respectively. In each problem, we consider the $\epsilon$-approximate solution with $\epsilon=0.1$ and $0.05$, respectively.

For each problem, we run 20 times and record the average results. Over 20 runs, the average values of the total number of iterations $N$ and the batch sizes $m$, as well as the smoothing parameters $\mu$ of the MSNS method are listed in Table \ref{tab:1}.

\begin{table}[h]
\caption{The average iteration limits, batch sizes and smoothing parameters over 20 runs}\label{tab:1}
\begin{center}
\begin{tabular}{llllll}
\hline\noalign{\smallskip}
$\epsilon$            & $NS$                   & $n$  & $N$ &$m$ &$\mu$\\
\noalign{\smallskip}\hline\noalign{\smallskip}
\multirow{4}{*}{0.1}  & \multirow{2}{*}{10000} & 500  &5275 &329 &0.0496 \\
                      &                        & 1000 &10587 &489&0.0299 \\
                      & \multirow{2}{*}{20000} & 500  &5292  &340 &0.0497 \\
                      &                        & 1000 &10513 &502&0.0316 \\
\noalign{\smallskip}\hline\noalign{\smallskip}
\multirow{4}{*}{0.05} & \multirow{2}{*}{10000} & 500  &21004  &675 &0.0251\\
                      &                        & 1000 &41851 &969 &0.0171\\
                      & \multirow{2}{*}{20000} & 500  &20906 &648 &0.0249 \\
                      &                        & 1000 &41679 &915 &0.0178\\
\noalign{\smallskip}\hline\noalign{\smallskip}
\end{tabular}
\end{center}
\end{table}

\begin{figure}[h]
\subfigure[ $NS=10000$, $n=500$, $\epsilon=0.1$]
{
\includegraphics[height=4.3cm]{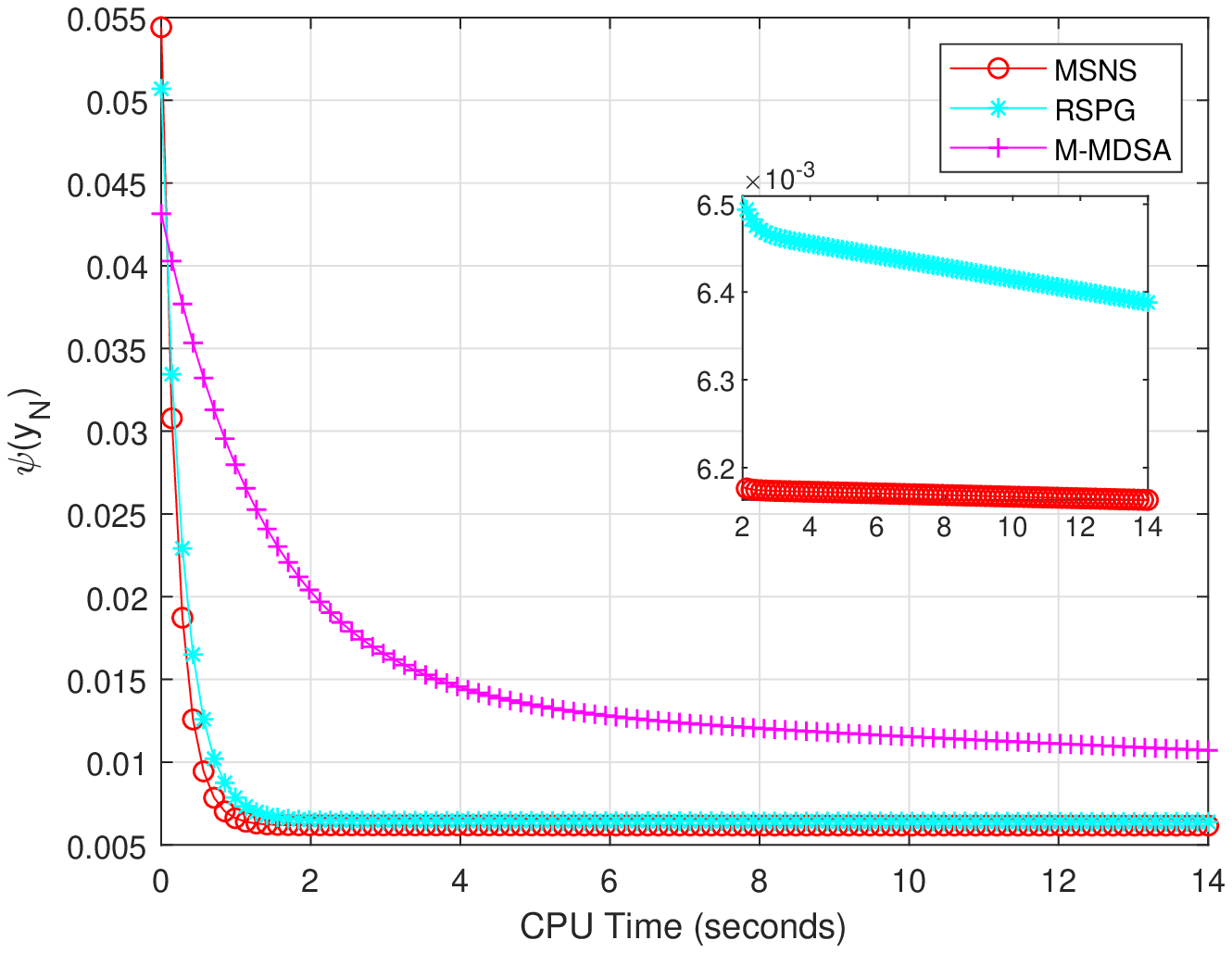}
}\quad
\subfigure[ $NS=20000$, $n=500$, $\epsilon=0.1$]
{
\includegraphics[height=4.3cm]{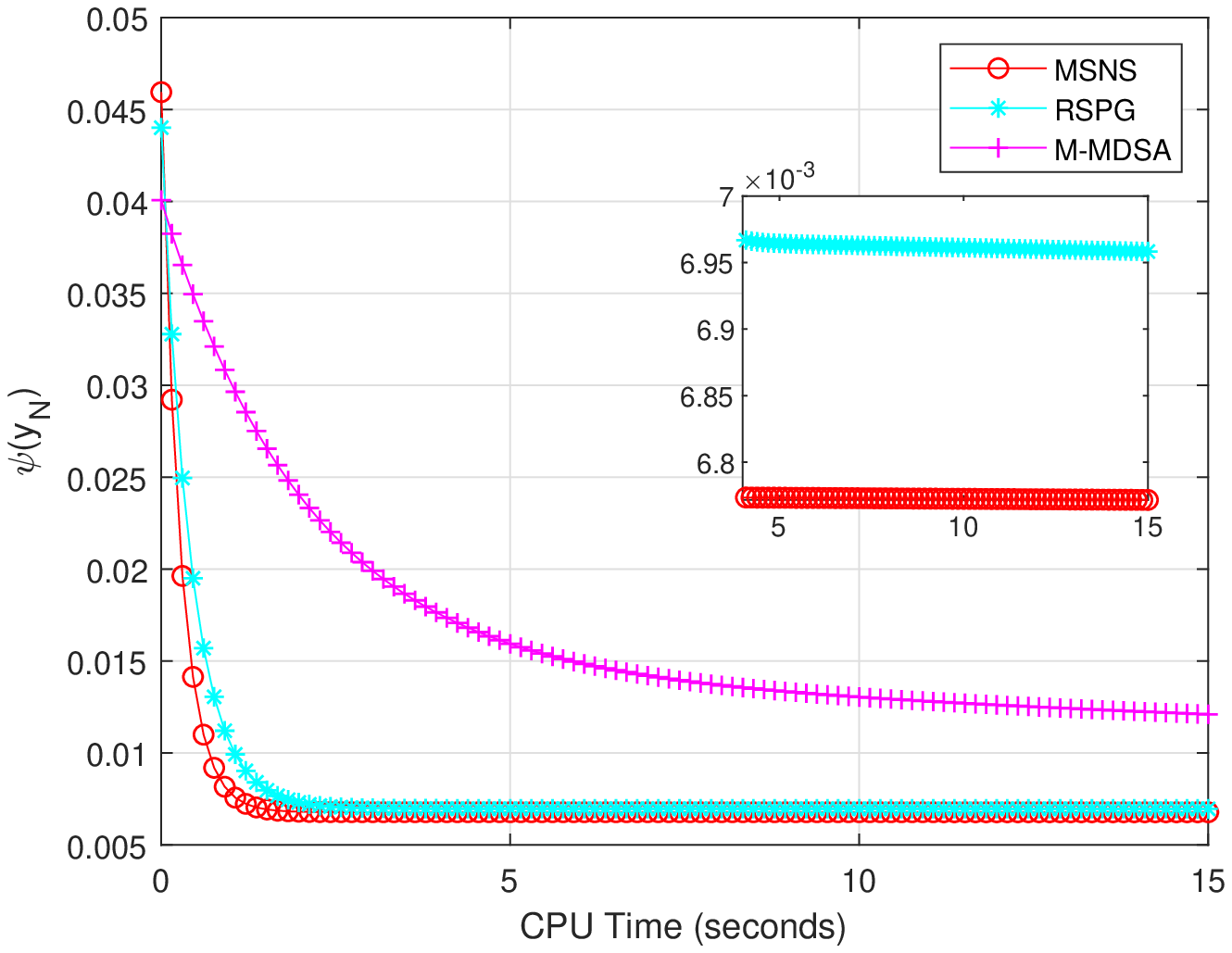}
}
\subfigure[ $NS=10000$, $n=1000$, $\epsilon=0.1$]
{
\includegraphics[height=4.3cm]{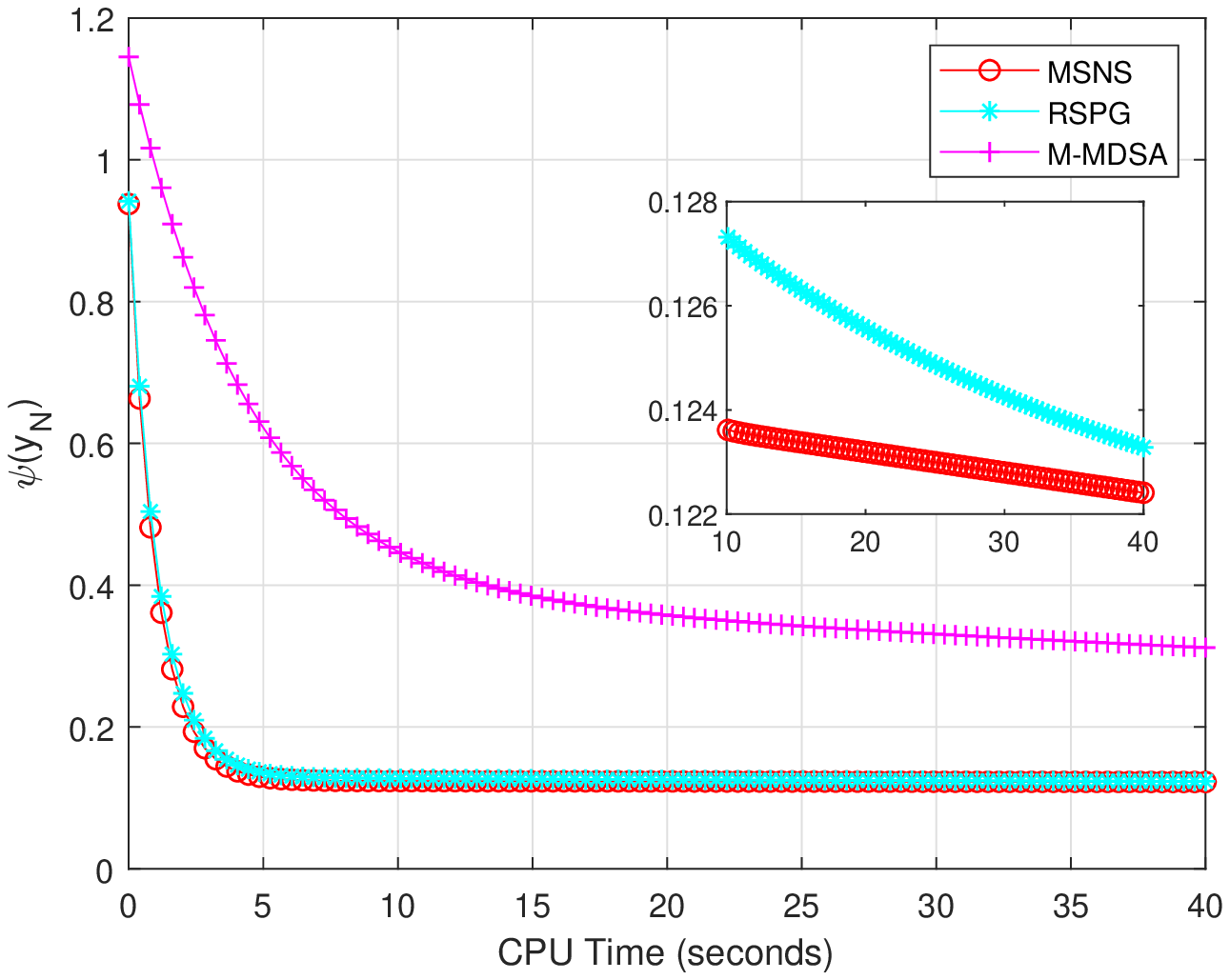}
}
\quad
\subfigure[$NS=20000$, $n=1000$, $\epsilon=0.1$ ]
{
\includegraphics[height=4.3cm]{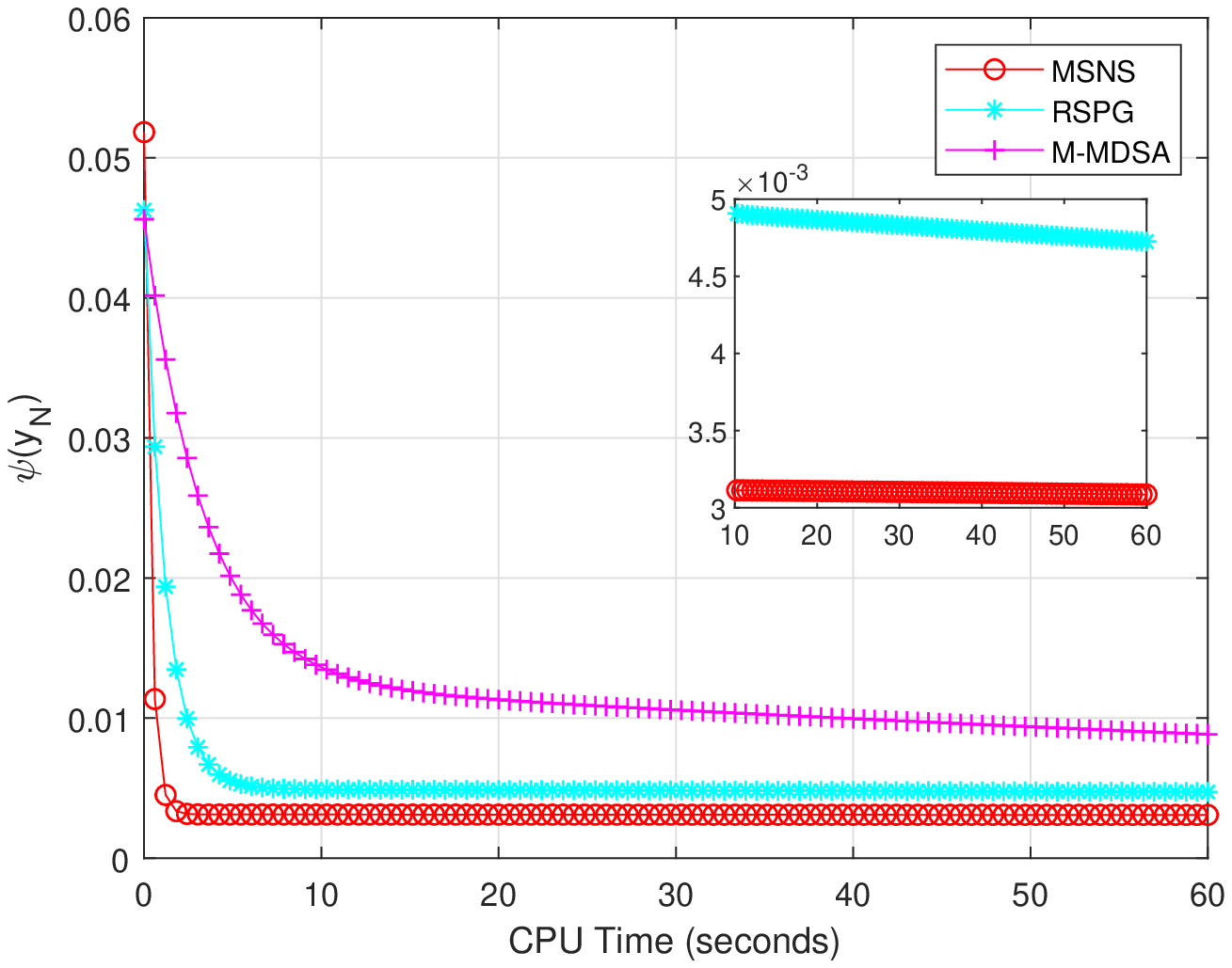}
}
\caption{Average objective values corresponding to the training data v.s. CPU time of 20 runs for different values of sample size $NS$, dimension $n$, when $\epsilon=0.1$.}\label{fig:1}
\end{figure}

\begin{figure}[h]
\subfigure[ $NS=10000$, $n=500$, $\epsilon=0.05$]
{
\includegraphics[height=4.3cm]{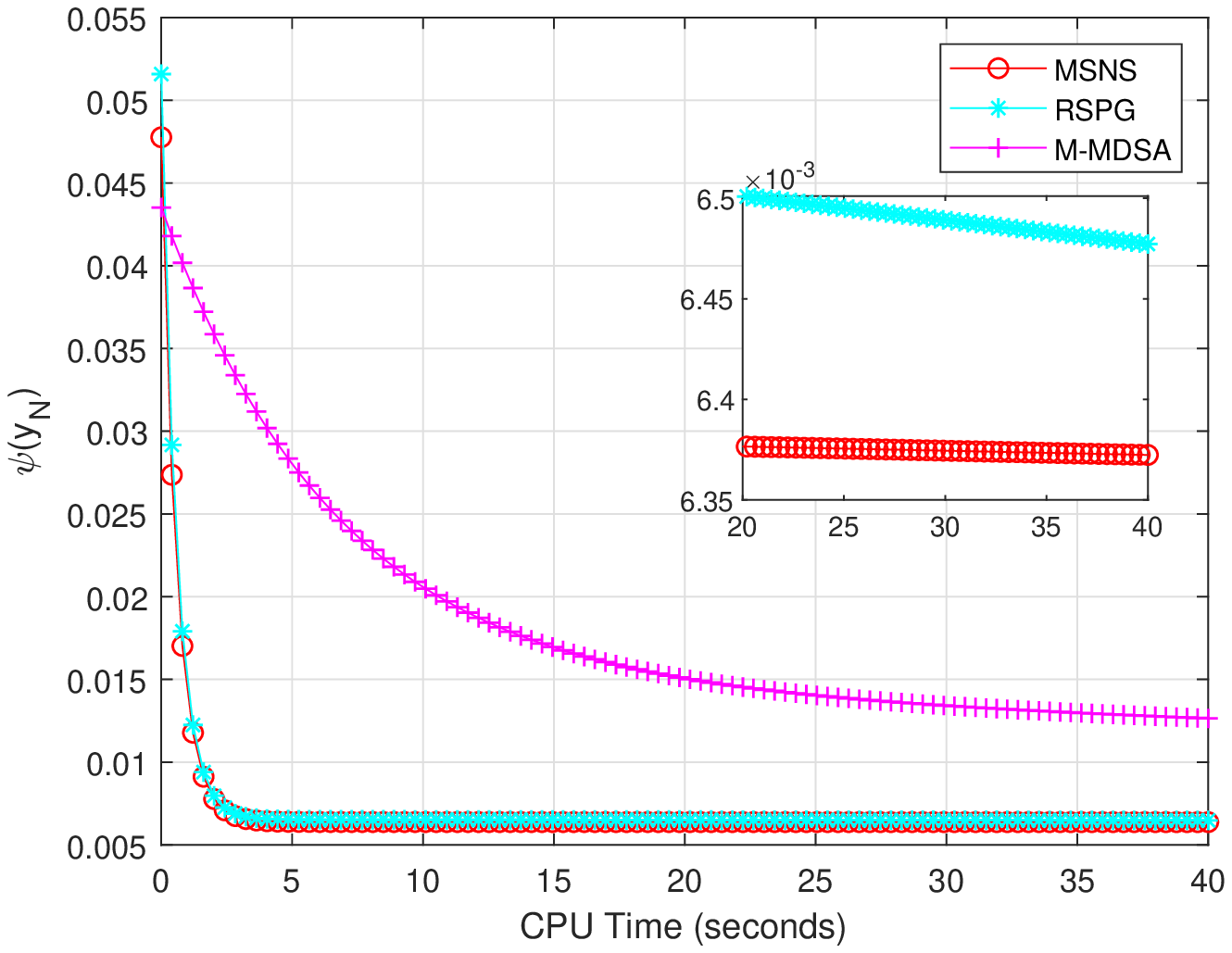}
}\quad
\subfigure[ $NS=20000$, $n=500$, $\epsilon=0.05$]
{
\includegraphics[height=4.3cm]{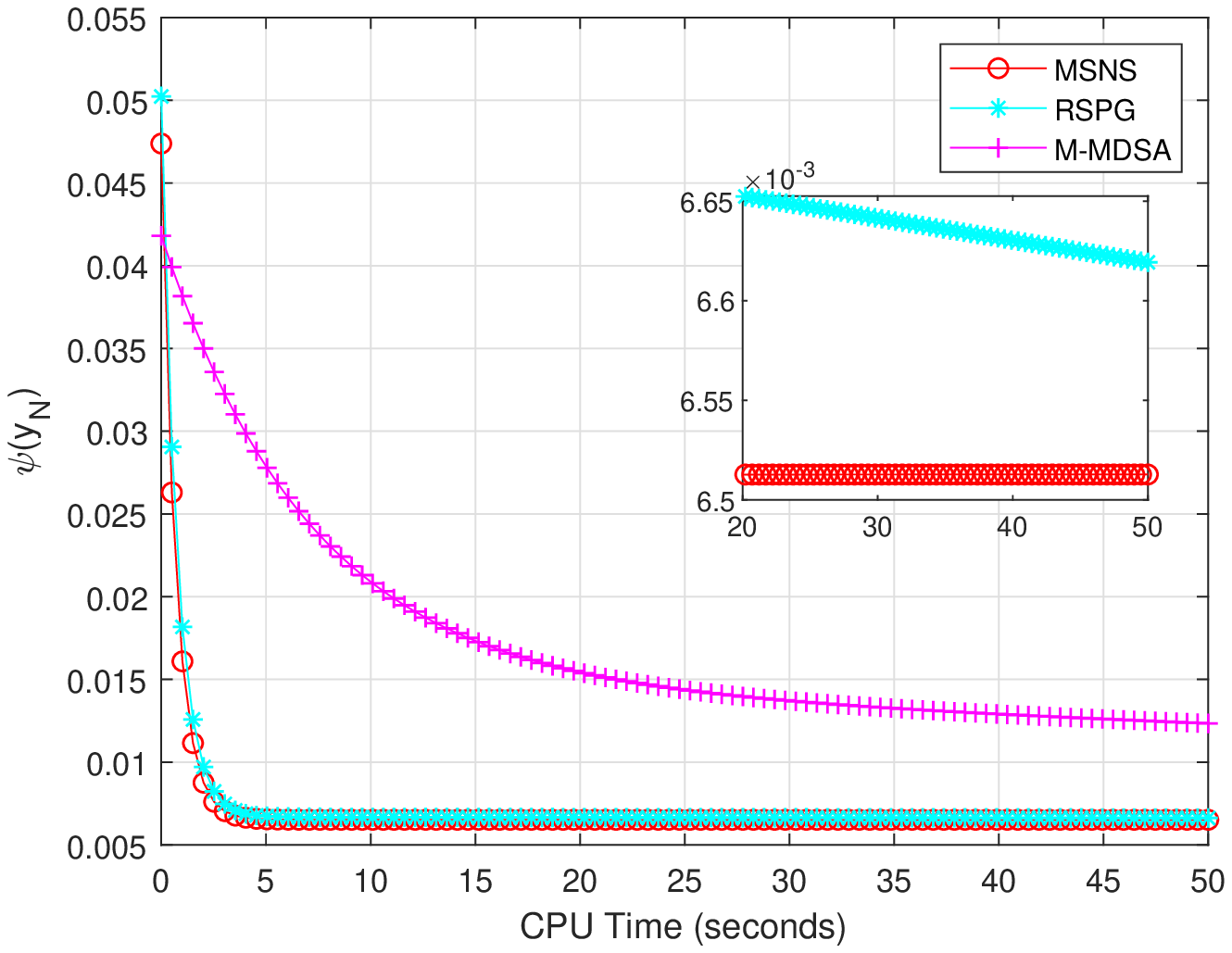}
}
\subfigure[ $NS=10000$, $n=1000$, $\epsilon=0.05$]
{
\includegraphics[height=4.3cm]{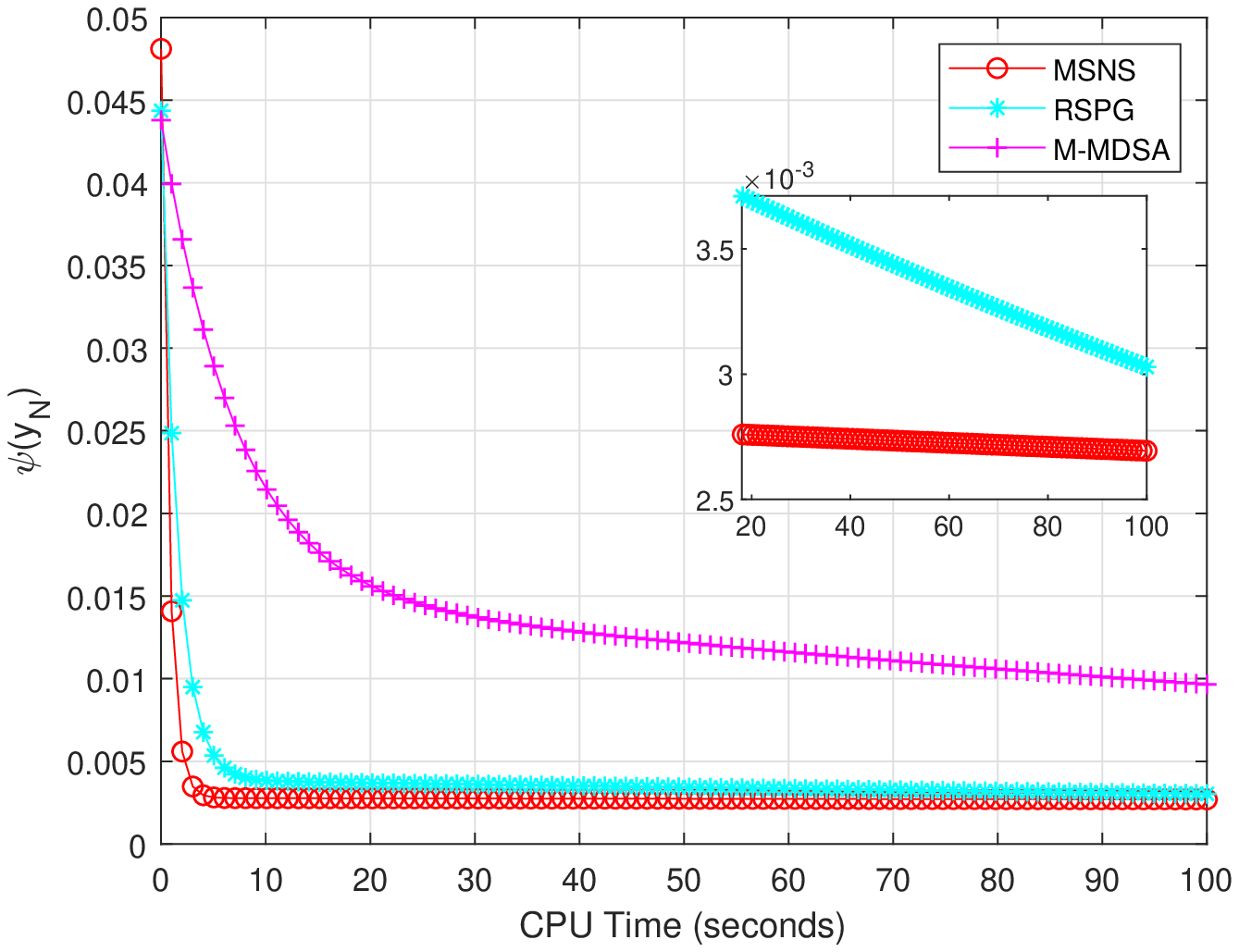}
}
\quad
\subfigure[$NS=20000$, $n=1000$, $\epsilon=0.05$ ]
{
\includegraphics[height=4.3cm]{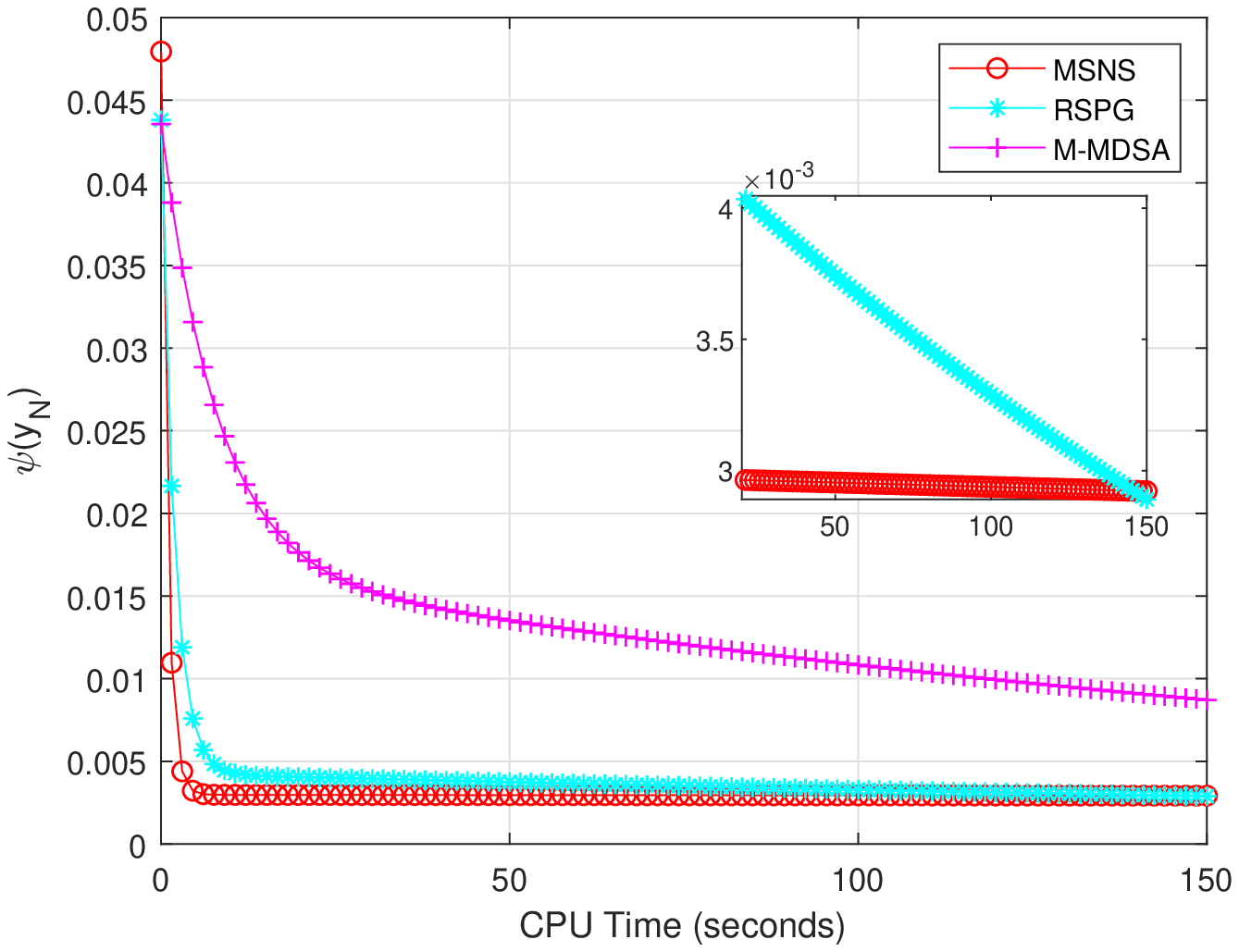}
}
\caption{Average objective values corresponding to the training data v.s. CPU time of 20 runs for different values of sample size $NS$, dimension $n$, when $\epsilon=0.05$.} \label{fig:2}
\end{figure}

We draw the curves of the average objective values corresponding to the training data v.s. the CPU time for the MSNS, RSPG and M-MDSA methods. We do not include the curves of the 2-RSPG and 2-RSPG-V methods in Fig. \ref{fig:1} and Fig. \ref{fig:2}, because they have several candidate outputs in optimization phase. Due to scalability, the curves of the MSNS and the RSPG methods are close at the end. We also provide small graph for each subfigure to see clear the differences of the two methods. We can see that the MSNS method provides the computed solution with the smallest objective values corresponding to the training data.

To see stability, we show in Tables \ref{tab:2} and \ref{tab:3} the mean (Mean) and the variance (Var), over 20 runs, of the objective values (Obj) corresponding to the testing data at the computed solution $\widehat{x}$ by a certain method. The Obj is the empirical mean of the stochastic objective value $F(\widehat{x}, \xi)+H(\widehat{x},\xi)$. The empirical mean is taken over a large testing data consisting of data samples $K = 50000$ as done in \cite{Nemirovski}. From Tables \ref{tab:2} and \ref{tab:3}, the average objective values corresponding to the testing data, the average accuracy,  and the average CPU time of the MSNS method are significantly better than that of the other methods in almost all cases. Moreover, the MSNS method provides the results of the objective values with small variance.

\begin{table}[h]
\caption{The mean and variance of the objective values, the average accuracy and CPU time when $\epsilon=0.1$  of 20 runs}
\label{tab:2}
\begin{center}
\begin{tabular}{lllllll}
\hline\noalign{\smallskip}
\multirow{2}{*}{$n$} & \multirow{2}{*}{$NS$}  & \multirow{2}{*}{ALG.} &{Obj} & {Obj} &Acc  & CPU   \\
                 &      &   & Mean  &Var   & (Avg.) & (Avg.) \\ \hline\noalign{\smallskip}
\multirow{10}{*}{500} &\multirow{5}{*}{10000}
&MSNS&\textbf{0.3071}&\textbf{2.90e-06}&\textbf{0.9817}&\textbf{44.2968}\\
 & & RSPG&0.3093&3.92e-06&0.9815&386.0023\\
 & & M-MDSA&0.3402&5.01e-05&0.9785&55.4375\\
 & & 2-RSPG&0.3092&2.05e-06&0.9815&334.7219\\
 & & 2-RSPG-V&0.3089&3.35e-06&0.9815&336.8773\\
 & \multirow{5}{*}{20000}
& MSNS&\textbf{0.3001}&\textbf{6.95e-07}&\textbf{0.9818}&\textbf{45.4281}\\
 & & RSPG&0.3031&7.95e-07&0.9813&386.8484\\
 & & M-MDSA&0.3357&2.00e-05&0.9783&57.9516\\
 & &2-RSPG&0.3028&1.05e-03&0.9814&424.7093\\
 & &2-RSPG-V&0.3027&9.57e-07&0.9816&437.3919\\ \hline\noalign{\smallskip}
\multirow{10}{*}{1000}&\multirow{5}{*}{10000}
&MSNS&\textbf{0.1715}&1.64e-06&\textbf{0.9978}&\textbf{102.5140}\\
 & & RSPG&0.1741&1.61e-06&0.9975&1335.9727\\
 & & M-MDSA&0.2351&1.21e-04&0.9924&103.5156\\
 & & 2-RSPG&0.1741&1.71e-06&0.9976&1669.5687 \\
 & & 2-RSPG-V&0.1739&\textbf{1.51e-06}&0.9976&1705.7171\\
 & \multirow{5}{*}{20000}
& MSNS&\textbf{0.1603}&\textbf{3.19e-07}&\textbf{0.9979}&\textbf{105.4656}\\
 & & RSPG&0.1646&6.66e-06&0.9976&1687.0039\\
 & & M-MDSA&0.2293&2.71e-05&0.9924&105.9687\\
 & &2-RSPG&0.1644&7.01e-07&0.9977 &1768.0484\\
 & &2-RSPG-V&0.1643&7.78e-07&0.9977 &1824.6265\\\hline\noalign{\smallskip}
\end{tabular}
\end{center}
\end{table}

\begin{table}[h]
\caption{The mean and variance of the objective values, the average accuracy and CPU time when $\epsilon=0.05$  of 20 runs}
\label{tab:3}
\begin{center}
\begin{tabular}{lllllll}
\hline\noalign{\smallskip}
\multirow{2}{*}{$n$}&\multirow{2}{*}{$NS$}  & \multirow{2}{*}{ALG.} &{Obj} & Obj &Acc  &CPU   \\
                 &      &   & Mean  &Var   & (Avg.) & (Avg.) \\ \hline\noalign{\smallskip}
\multirow{10}{*}{500} &\multirow{5}{*}{10000}
&MSNS&\textbf{0.3061}&1.71e-06&\textbf{0.9813}&\textbf{266.1671}\\
 & &RSPG&0.3076&1.85e-06&0.9811 &1951.5351\\
 & &M-MDSA&0.3411&4.76e-05&0.9781&391.4593\\
 & &2-RSPG&0.3074&\textbf{1.61e-06}&0.9812&1848.8625\\
 & &2-RSPG-V&0.3072&1.88e-06&0.9812&1835.7265\\
 &\multirow{5}{*}{20000}
& MSNS&\textbf{0.2995}&1.47e-06&\textbf{0.9823}&\textbf{273.3210}\\
 & & RSPG&0.3008&1.68e-06&0.9820&2693.3601\\
 & & M-MDSA&0.3324&1.91e-05&0.9791&399.0406\\
 & &2-RSPG&0.3008&\textbf{1.36e-06}&0.9821&1915.3305\\
 & &2-RSPG-V&0.3007&1.63e-06&0.9821&2462.5796\\ \hline\noalign{\smallskip}
\multirow{10}{*}{1000}&\multirow{5}{*}{10000}
&MSNS&\textbf{0.1704}&\textbf{1.42e-06}&\textbf{0.9977}&\textbf{514.4062}\\
 & &RSPG&0.1721 &1.48e-06&0.9975&10347.2710\\
 & &M-MDSA&0.2364&5.96e-05&0.9923&529.8718\\
 & &2-RSPG&0.1719&1.95e-06&0.9976 &8604.6609\\
 & &2-RSPG-V&0.1717&1.78e-06&0.9976 &9281.1523\\
 &\multirow{5}{*}{20000}
& MSNS&\textbf{0.1595}&\textbf{4.61e-07}&\textbf{0.9979} &\textbf{534.6511}\\
 & & RSPG&0.1614&6.66e-07&0.9978 &7236.1585\\
 & & M-MDSA&0.2283&3.21e-05&0.9924 &541.6226\\
 & &2-RSPG&0.1616&6.88e-07&0.9978 &9701.7789\\
 & &2-RSPG-V&0.1614&5.32e-07&0.9978&9028.8843\\ \hline\noalign{\smallskip}
\end{tabular}
\end{center}
\end{table}

\subsection{Real datasets}
We do numerical experiments on four real datasets described below for our experiments.

$\bullet$ Wisconsin breast cancer dataset from the UCI repository (699 patterns) can be downloaded from the web\footnotemark[1]\footnotetext[1]{https://archive.ics.uci.edu/ml/datasets/Breast+Cancer+Wisconsin+(Diagnostic)}. Features are computed from a digitized image of a fine needle aspirate (FNA) of a breast mass. They describe characteristics of the cell nuclei present in the image.

$\bullet$ Statlog (Australian Credit Approval) dataset also comes from the UCI repository, downloaded from the web\footnotemark[2]\footnotetext[2]{https://archive.ics.uci.edu/ml/datasets/Statlog+\%28Australian+Credit+Approval\%29}. This file concerns credit card applications. All attribute names and values have been changed to meaningless symbols to protect confidentiality of the data. This dataset is interesting because there is a good mix of attributes -- continuous, nominal with small numbers of values, and nominal with larger numbers of values.

$\bullet$ Credit Approval dataset also comes from the UCI repository\footnotemark[3]\footnotetext[3]{https://archive.ics.uci.edu/ml/datasets/Credit+Approval}. This file concerns credit card applications. All attribute names and values have been changed to meaningless symbols to protect confidentiality of the data. This dataset is interesting because there is a good mix of attributes, continuous, nominal with small numbers of values, and nominal with larger numbers of values.

$\bullet$ Ecoli dataset is also refer to the protein localization sites dataset, downloaded from the web\footnotemark[4]\footnotetext[4]{https://archive.ics.uci.edu/ml/datasets/Ecoli}. The dataset describes the problem of classifying Ecoli proteins using their amino acid sequences in their cell localization sites. That is, predicting how a protein will bind to a cell based on the chemical composition of the protein before it is folded. We analyzed Ecoli dataset in its 2-class versions, i.e., Ecoli(B). Ecoli(B) from the first 4 proteins and the remaining ones. The details of the described datasets are resumed in Table \ref{tab:4}.

\begin{table}[h]
\caption{Details of the datasets}
\label{tab:4}
\begin{center}
\begin{tabular}{llll}
\hline
Dataset  & Classes &Sample size &Dimension        \\ \hline\noalign{\smallskip}
Wisconsin breast cancer  & 2       & 699 & 10  \\
Statlog  & 2       & 690 & 14  \\
Credit Approval  & 2       & 690 & 15  \\
Ecoli(B)  & 2       & 366 & 343  \\ \hline\noalign{\smallskip}
\end{tabular}
\end{center}
\end{table}

We choose the optimal values of $\lambda_1$ and $t$ via 3-fold cross-validation (CV) using 20 random runs, which are determined by varying them on the grid $\{10^{-2},10^{-1},2^{-2},2^{-1},1\}$ and the values with the best average accuracy are chosen for each of the MSNS, RSPG, M-MDSA, 2-RSPG, and 2-RSPG-V methods. For real datasets, we also compare with the classical SVM model \cite{chang} for which the cost parameter $c$ and parameter $g$ in kernel function were determined by varying them on the grid $\{10^{-4},10^{-2},10^{-1},1,10\}$ recommended by \cite{Melacci} on page 1168.

We find that the average CPU time is related to the parameters $\lambda_1$ and $t$ for each of the MSNS, RSPG, M-MDSA, 2-RSPG, and 2-RSPG-V methods. More specifically, the CPU time decreases as $\lambda_1$ decreases, and also decreases as $t$ decreases for each of the above method. We take the MSNS method as an example to illustrate the reason for this phenomenon. When $\lambda_1$ is fixed, along with the decreasing of $t$, $D=t/2$  decreases and consequently the maximum iteration number $N$ decreases according to \eqref{complexity}. The batch size $m$ also decreases with $N$ by \eqref{eq3.34}. Similarly, when $t$ is fixed, $L_f$ decreases with $\lambda_1$, resulting in a decrease in the maximum number of iterations $N$ by \eqref{complexity}. The batch size $m$ also decreases with $N$ by \eqref{eq3.34}. We record the average CPU time of the MSNS method, along with different $\lambda_1$ and $t$  for two datasets: Wisconsin breast cancer and Ecoli(B) as examples in Tables \ref{tab:5} and \ref{tab:6}, respectively. The relationship of the CPU time and the values of parameters $\lambda_1$ and $t$ can be seen clearly from the two tables.

\begin{table}[h]
\caption{CPU time of the MSNS method corresponding to different values of $\lambda_1$, $t$, on the dataset -- Wisconsin breast cancer}
\label{tab:5}
\begin{center}
\begin{tabular}{llllll}
\hline
$t\backslash\lambda_1$&$10^{-2}$&$10^{-1}$&$2^{-2}$&$2^{-1}$&1\\ \hline\noalign{\smallskip}
 $10^{-2}$&0.0074  &0.0038  &0.0040  &0.0043  &0.0060  \\
 $10^{-1}$&0.0246  &0.0259  &0.0310  &0.0366  &0.0453  \\
 $2^{-2}$&0.0653  &0.0661  &0.0784  &0.0942  &0.1231  \\
 $2^{-1}$&0.1258  &0.1433  &0.1631  &0.2079  &0.2725  \\
 1&0.2676  &0.3013  &0.3574  &0.4541  &0.6094  \\ \hline\noalign{\smallskip}
\end{tabular}
\end{center}
\end{table}

\begin{table}[h]
\caption{CPU time of the MSNS method corresponding to different values of $\lambda_1$, $t$, on the dataset -- Ecoli(B)}
\label{tab:6}
\begin{center}
\begin{tabular}{llllll}
\hline\noalign{\smallskip}
$t\backslash\lambda_1$&$10^{-2}$&$10^{-1}$&$2^{-2}$&$2^{-1}$&1\\ \hline\noalign{\smallskip}
 $10^{-2}$ &0.3350 &0.3659 &0.3885 &0.4317&0.5288\\
 $10^{-1}$&9.8397 &10.2583 &10.9856 &12.2976&15.0089\\
 $2^{-2}$&39.5078 &41.1832 &44.0774 &48.2053 &58.5466\\
 $2^{-1}$&103.0333  &114.9313 &127.1659 &134.4300 &163.9350 \\
 1 &211.7660 &285.0719 &322.8352 &399.9200&471.5500\\ \hline\noalign{\smallskip}
\end{tabular}
\end{center}
\end{table}

We record in Table \ref{tab:7} the average optimal values of parameters $(t,\lambda_1)$ or $(c,g)$, together with the corresponding average accuracy and the average CPU time in seconds. We can see that our MSNS method has the best average accuracy in all the datasets. The CPU time of the MSNS method is not the shortest, but it is only a little bit longer than the shortest one and hence is acceptable.

\begin{table}[h]
\caption{Accuracy, the values of $(t,\lambda_1)$ or $(c,g)$ and CPU time determined by 3-fold CV}
\label{tab:7}
\begin{center}
\begin{tabular}{llllllll}
\hline\noalign{\smallskip}
Dataset
& ALG.     & $t$ &$\lambda_1$& $c$ &$g$ & Acc & CPU      \\ \hline
\multirow{6}{*}{\begin{tabular}[c]{@{}c@{}}Wisconsin\\ breast cancer\end{tabular}}
& MSNS     &$10^{-1}$ &$10^{-2}$&-&- &\textbf{0.9686}   & 0.0097  \\
& RSPG     &$10^{-1}$ &$10^{-1}$&-&- &0.9664  &0.0597  \\                                                                                           & M-MDSA   &$10^{-1}$ &$10^{-2}$&-&- &0.9617  &0.0111 \\
& 2-RSPG   &$10^{-1}$ &$10^{-2}$&-&- &0.9672  &0.0674 \\
& 2-RSPG-V &$10^{-1}$ &$10^{-2}$&-&- &0.9678  &0.0566\\
& SVM     &-&- &10&$10^{-2}$ &0.9614   &\textbf{0.0046}   \\ \hline\noalign{\smallskip}
\multirow{6}{*}{\begin{tabular}[c]{@{}c@{}}Statlog\end{tabular}}                                                                  & MSNS     &$2^{-2}$&$10^{-2}$&-&- & \textbf{0.8670} &\textbf{0.0163}\\
& RSPG     &$2^{-2}$&$10^{-2}$&-&-&0.8631  & 0.0487\\
& M-MDSA   & 1&$10^{-1}$&-&-& 0.8569 & 0.1602\\
& 2-RSPG   &$2^{-2}$&$10^{-2}$ &-&-&0.8654 & 0.0765\\
& 2-RSPG-V &$2^{-2}$&$10^{-2}$&-&- &0.8656   &0.0580 \\
& SVM      &-&-& 10 & $10^{-4}$& 0.8609&0.0254   \\\hline\noalign{\smallskip}
\multirow{6}{*}{\begin{tabular}[c]{@{}c@{}}Credit Approval\end{tabular}}                                                                  & MSNS     &$2^{-1}$&1&-&- &\textbf{0.8595}   &0.2636 \\
& RSPG     &$2^{-2}$&1 &-&-&0.8589 &0.4914\\
& M-MDSA   &1&$2^{-1}$&-&-&0.8562&0.4694\\
& 2-RSPG   &1 &1 &-&-&0.8591 &2.9523\\
& 2-RSPG-V &$2^{-2}$&1&-&- &0.8590&0.5832\\
& SVM      &-&-& 1 &$10^{-2}$ &0.8591 &\textbf{0.0155}   \\\hline\noalign{\smallskip}
\multirow{6}{*}{\begin{tabular}[c]{@{}c@{}}Ecoli(B)\end{tabular}}                                                                  & MSNS     &$10^{-2}$&$10^{-1}$&-&-&\textbf{0.8720}   &0.3659  \\
& RSPG     &$10^{-2}$&$2^{-1}$&-&-&0.8501&0.2826  \\
& M-MDSA   &$10^{-1}$&$10^{-1}$&-&-&0.8423  &10.2845 \\
& 2-RSPG   &1&$10^{-2}$&-&-&0.8540 &160.8485\\
& 2-RSPG-V &1&$10^{-2}$&-&-&0.8542  &169.9453 \\
& SVM   &-&-& 10& $10^{-4}$&0.7107  &\textbf{0.0767}  \\ \hline\noalign{\smallskip}
\end{tabular}
\end{center}
\end{table}

\section{Concluding remarks}
In this paper, we propose a mini-batch stochastic Nesterov's smoothing (MSNS) method for solving a class of constrained convex nonsmooth composite optimization problems with noisy zero-order and first-order information. We show the convergence of the MSNS method, together with its optimal iteration complexity. Numerical experiments on a support vector machine (SVM) model using both synthetic datasets and real datasets, demonstrate the effectiveness and efficiency of the proposed MSNS method, compared with several state-of-the-art methods.

\end{document}